\documentclass[12pt,leqno]{amsart}
\usepackage{amsfonts}
\usepackage{amssymb}
\usepackage{amsthm}
\usepackage{enumitem}
\usepackage{hyperref}
\usepackage[english]{babel}

\def\dist{{\rm dist}}
\def\mod{{\rm mod}}

\def\ds{\displaystyle}

\newtheorem{thm}{Theorem}
\newtheorem{lem}{Lemma}

\newtheorem{ex}{Example}

\newtheorem{prop}{Proposition}

\newtheorem{defn}{Definition}
\newtheorem{cor}{Corollary}

\begin{document}
	\title{Contraction map sets with an external factor and weakly fixed points}
	\author{Vasil Zhelinski}
	\email{vasil\_zhelinski@uni-plovdiv.bg}
	
	\maketitle
	\section{\normalsize Abstract}
	This research aims to introduce a new property, called the $CD$ property as a less restrictive and more user-friendly variant of the $UC$ one.
	This research then aims, by using the idea of noncyclic maps, to reveal a new method for generalizing Banach's fixed point theorem by introducing the concept of an external factor.
	Another goal of this research is to show a new way of obtaining results for $p$-cyclic maps by using the results of noncyclic maps. Also this research shows a relationship between fixed points and best proximity points similar to that between fixed points and coupled fixed points as shown in \cite{PETROSEL}.
	\section{\normalsize Introduction}
	The study of contraction maps begins with Banach's fixed point theorem \cite{BANACH}, which turns to be an effective tool both in applied and pure mathematics. This theorem refers self maps $T:X\to X$, defined on a complete metric spaces $(X,\rho)$ and satisfying the contractive $\rho(Tx,Ty)\leq \lambda \rho(x,y)$ for every $x,y\in X$ and some $\lambda\in [0,1)$.
	
	Over time an enormous number of generalizations have been made. We will try to classify them in few types.  One kind is to alter the underlying space by considering partial metric space first used in \cite{PMS-INTRODUCTED} and studied recently in \cite{PMS-3,PMS-2,PMS-1}, partially ordered metric spaces used first for coupled fixed points in \cite{Ran} and some results are \cite{Guo,PORDEREDMS,Lakshmikantham}, $b$-metric space at first used in \cite{B-MS-INTRODUCED} and considered recently in \cite{B-MS-1,B-MS-2,B-MS-3}, Banach spaces \cite{BPP,UCBS1,ZLATANOV,ZLATANOV2}, modular function spaces used first in \cite{KKR} and later researched in \cite{Ilchev20172873,Khamsi-Kozlowski,Kozlowski,ZLATANOV-m}, etc.
	
	Different approach is the changing of the type of contraction conditions e Kannan \cite{KANAN}, Chatterjea \cite{CHATTERJEA}, Hardy–Rogers \cite{HARDY-ROGERS}, Meir Keeler \cite{MK} maps, etc.
	
	Another approach to a generalization of Banach's theorem is to consider maps of two variables. In this context, the need for a suitable generalization of the concept of a fixed point occurs. This is the idea for coupled fixed points of a map of two variables $F$, introduced in \cite{COUPLEDFP}, defined as $(x,y)$ such $F(x,y)=x$, $F(y,x)=y$. The fixed points for multi variable maps are a generalization of coupled fixed points for maps of $n$ variables.
	
	Important results about the connection between coupled fixed points and fixed points
	are obtained in \cite{PETROSEL}
	
	A recent approach in the generalization of to replace the self map with a cyclic $T(A)\subseteq B$, $T(B)\subseteq A$ introduced in \cite{CYCLIC} or with a noncyclic ones \cite{NONCYCLIC}, $T(A)\subseteq A$, $T(B)\subseteq B$. As far as if the maps $T$ is a cyclic one and there may hold ${\rm dist}(A,B)=\inf\{\rho(a,b):a\in A,b\in B\}>0$ it is proposed in  \cite{BPP} to search for a point $x\in A\cup B$ so that to minimize the function $\rho(x,Tx)$. If there holds $\rho(x,Tx)=0$ then $x$ is a fixed point of the map $T$. A sufficient condition is found in \cite{BPP} so that cyclic maps to satisfy $\rho(x,Tx)=\min\{\rho(x,Tx):x\in A\cup B\}={\rm dist}(A,B)$.
	The idea to consider cyclic maps have been generalized in defining $p$-cyclic maps introduced in \cite{PCYCLIC}. Recent results on the topic are \cite{P-CYCLIC-1,P-CYCLIC-2,P-CYCLIC-4,P-CYCLIC-5,P-CYCLIC-3}.
	Let $\{A_i\}_{i=1}^p$, where $p>1$, be nonempty subsets of a metric space $(X, \rho)$. We use the convention $A_{p+i} = A_i$ for every $i \in \{1, 2, 3, \dots\}$
	The map $\displaystyle T :\bigcup_{i=1}^p A_i\to\bigcup_{i=1}^p A_i$
	is called a $p$-cyclic map if $T(A_i)\subseteq A_{i+1}$ for every $i = 1, 2,\dots p$. We easily can see that for $p=2$ the definition of a $2$-cyclic map matches that of a cyclic map. Point $a\in A_i$ is called a best proximity point of the p-cyclic map $T$ in $A_i$ if $\rho(a,Ta)= \dist(A_i,A_{i+1})$. An interesting fact about $p$-cyclic maps is that in the start, the theorems carry with themselves the restriction that $\dist(A_1,A_2)=\dist(A_2,A_3)=\dots=\dist(A_{p-1},A_p)=\dist(A_{p},A_1)$ an example of this is \cite{PCYCLIC}, but later this restriction is overcame, this can be seen in \cite{UCBS1,ZLATANOV,ZLATANOV2}.
	
	There are many recent results on the topic of noncyclic maps \cite{NONCYCLIC3,NONCYCLIC1,NONCYCLIC2,NONCYCLIC4,NONCYCLIC5}.
	To finish the listing of the possible generalization let us say that we can consider all  mentioned types simultaneously, an example of combining more than one type is \cite{Koleva}.  
	
	Most of the results on existence and uniqueness of best proximity points for cyclic or $p$-cyclic maps rely on the good geometry of the underlying Banach space \cite{BPP,UCBS1,ZLATANOV,ZLATANOV2}. The key lemmas from \cite{BPP} that hold in uniformly convex Banach space a used in the proofs of existence and uniqueness of best proximity points. Later it is suggested to replace the key lemmas from \cite{BPP} by the properties $UC$ defined in \cite{UC}, $UC^{*}$ introduced in \cite{UCSTAR}, the $P$-property, used in \cite{P-P1,P-P2}.
	The relationship between $UC$ and $UC^{*}$ and the uniform convexity of a Banach space, is explored in \cite{UC-UCZ}.
	Recent results using the property $UC$ can be found in \cite{UC-5,UC-4,UC-3,UC-2,UC-1}.
	
	This article aims to introduce a new property, called the $CD$ property as a less restrictive and more user-friendly variant of the $UC$ one.
	
	This research then aims, by using the idea of noncyclic maps, to reveal a new method for generalizing Banach's fixed point theorem by introducing the concept of an external factor.
	Another goal of this research is to show a new way of obtaining results for $p$-cyclic maps by using the results of noncyclic maps. Also this research shows a relationship between fixed points and best proximity points similar to that between fixed points and coupled fixed points as shown in \cite{PETROSEL}.
	
	\section{\normalsize Preliminaries}
	
	In this research we will use the notations $\mathbb{N}$ for the set of natural numbers, $\mathbb{N}_0$ for $\mathbb{N}\cup\{0\}$,
	$\mathbb{R}$ for the set of real numbers, $\mathbb{Z}$ for the set of integers,\\
	$\lfloor\cdot\rfloor$ for the floor function ($\lfloor x \rfloor=\max\{n\in \mathbb{Z}:n\leq x\}$), $\lceil\cdot\rceil$ for the ceiling function ($\lceil x \rceil=\min\{n\in \mathbb{Z}:n\geq x\}$), $(X,\rho)$ for a metric space $X$ with metric $\rho$, $(X,\|\cdot\|)$ for a normed space $X$ with norm $\|\cdot\|$. We will consider the
	metric generated by the norm. For the unit sphere in a
	normed space $(X,\|\cdot\|)$ we will use the notation $S_X$, ($S_X=\{x\in X: \|x\|=1\}$).
	
	Let $(X, \rho)$ be a metric space, $A$ and $B$ be subsets of $X$ we will denote the distance between $A$ and $B$ by $\displaystyle\dist(A,B)=\inf_{a\in A,b\in B}\rho(a,b)$.
	
	\begin{defn}[\cite{UC}]
		Let $(X,\rho)$ be a metric space, $A$ and $B$ be subsets of $X$. The ordered pair $(A, B)$ satisfies the property $UC$, if for every triple sequences $\{x_n\}_{n=0}^{\infty}\subseteq A$, $\{z_n\}_{n=0}^{\infty}\subseteq A$ and $\{y_n\}_{n=0}^{\infty}\subseteq B$ such that
		\begin{itemize}
			\item $\displaystyle \lim_{n\to \infty}\rho(x_n,y_n)=\dist(A,B)$
			\item $\displaystyle \lim_{n\to \infty}\rho(z_n,y_n)=\dist(A,B)$
		\end{itemize}equality
		there holds $\displaystyle \lim_{n\to \infty}\rho(x_n,z_n)=0$
	\end{defn}
	
	\begin{defn}[\cite{ZLATANOV}]
		Let $\{A_i\}_{i=1}^3$ be nonempty subsets of a Banach space  $(X, \|\cdot\|)$,  $\displaystyle T :\bigcup_{i=1}^3 A_i\to\bigcup_{i=1}^3 A_i$ be a cyclic map and there exists $k\in (0,1)$ so that for every $x_1\in A_1$, $x_2\in A_2$, $x_3\in A_3$ there holds
		\begin{equation}\label{bz}
			\begin{array}{rcl}
				W_1&=&\|Tx_1-Tx_2\|+\|Tx_2-Tx_3\|+\|Tx_3-Tx_1\|\\
				&\leq&k(\|x_1-x_2\|+\|x_2-x_3\|+\|x_3-x_1\|)+(1-\lambda)D
			\end{array}
		\end{equation}
		where $D= \dist(A_1,A_2)+ \dist(A_2,A_3)+ \dist(A_3,A_1)$.
		
		Then we say that $T$ is a 3-cyclic summing contraction
	\end{defn}
	
	\begin{defn}[\cite{ZLATANOV}]
		Let $(X,\|\cdot\|)$ be a Banach space. The function\\
		$\delta_{\|\cdot\|}:[0;2]\to[0;1]$, defined as\\
		$\ds \delta_{\|\cdot\|}(\varepsilon)=\inf\left\{1-\left\|\frac{x+y}{2}\right\|:x\in S_X, y\in S_X,\|x-y\|\geq \varepsilon\right\}$, is called modulus of convexity of $(X,\|\cdot\|)$. If $\delta_{\|\cdot\|}(0,2]\subseteq(0,1]$, then we say that $(X,\|\cdot\|)$ is a uniformly convex Banach space.
	\end{defn}
	
	\begin{ex}[\cite{UC}]\label{e0}
		The following are examples of a pairs of nonempty subsets $(A,B)$
		satisfying the property $UC$.
		\begin{enumerate}[label=\arabic*)]
			\item Every pair of nonempty subsets $(A, B)$ of a metric space $(X, \rho)$ such that $ \dist(A, B) = 0$.
			\item Every pair of nonempty subsets $(A, B)$ of a uniformly convex Banach space $(X,\|\cdot\|)$ such
			that $A$ is convex.
			\item Every pair of nonempty subsets $(A, B)$ of a strictly convex Banach space which A is
			convex and relatively compact and the closure of $B$ is weakly compact.
		\end{enumerate}
	\end{ex}
	
	Let us point out that the modulus of convexity $\delta$ is considered if the Banach space is at least two dimensional.
	A deep observation \cite{ERROR-1}, \cite{ERROR-2}, \cite{ERROR-3} allow us to get results about best proximity points in  $(\mathbb{R},|\cdot|)$. 
	
	\begin{cor}\label{c-0}
		If $A$ and $B$ are real intervals, then the ordered pair $(A,B)$ satisfies the $UC$ property.
	\end{cor}
	
	\begin{thm}[\cite{ZLATANOV}]\label{t-1}
		Let $\{A_i\}_{i=1}^3$ be closed, convex subsets of an uniformly convex Banach space $(X, \|\cdot\|)$,  $\displaystyle T :\bigcup_{i=1}^3 A_i\to\bigcup_{i=1}^3 A_i$ be a 3-cyclic summing contraction. Then there exist unique points $z_i\in A_i$ so that 
		\begin{itemize}
			\item for any $x\in A_i$ the sequence $\{T^{3n}x\}_{n=1}^{\infty}$ converges to $z_i$
			\item $\rho(z_1,z_2)= \dist(A_1,A_2)$, $\rho(z_2,z_3)= \dist(A_2,A_3)$ and\\
			$\rho(z_3,z_1)= \dist(A_3,A_1)$
			\item $Tz_1=z_2$, $Tz_2=z_3$, $Tz_3=z_1$
		\end{itemize}
		
	\end{thm}
	
	\begin{defn}\label{def4}
		Let $A$ and $C$ be some sets, let $T:A\times C \to A$ and  $H:A\times C \to C$, let $p=(x_0,u_0)\in A\times C$, $\{x_n\}_{n=0}^\infty$ and $\{u_n\}_{n=0}^\infty$ be defined as $x_{n+1}=T(x_n,u_n)$, $u_{n+1}=H(x_n,u_n)$
		
		Then we say that $\{x_n\}_{n=0}^\infty$ and $\{u_n\}_{n=0}^\infty$ are the iterated sequences of $(T,H)$ with an initial guess $p$.
		
	\end{defn}
	
	\section{\normalsize Main Results}
	
	Now we will introduce a property that is a more convenient variant of the $UC$ property, namely:
	
	\begin{defn}
		Let $(X,\rho)$ be a metric space, $A$ and $B$ be subsets of $X$, $C$ be an arbitrary set, $P\subseteq A\times B \times C^2$, $f_A:C\to \mathbb{R}$, $f_B:C\to \mathbb{R}$ and for every two sequences, $\{x_n\}_{n=1}^{\infty}\subseteq A$ and $\{y_n\}_{n=1}^{\infty}\subseteq B$ such that 
		\begin{itemize}
			\item $\displaystyle\lim_{k\to \infty}\sup_{n\geq k,\ m\geq k}\rho(x_n,y_m)= \dist(A,B)$
			\item there exist $\{u_n\}_{n=1}^{\infty}\subseteq C$ and $\{v_n\}_{n=1}^{\infty}\subseteq C$ so that $(x_n,y_m,u_n,v_m)\in P$ for every $n,m\in \mathbb{N}$, $\displaystyle \lim_{n\to \infty}f_A(u_n)=\inf_{c\in C}f_A(c)$ and $\displaystyle \lim_{n\to \infty}f_B(v_n)=\inf_{c\in C}f_B(c)$
		\end{itemize} there holds $\{x_n\}_{n=1}^{\infty}$ converges in $A$. Then we say that the ordered pair $(A,B)$ satisfies the property Convergent by Distance in $P$ on $(f_A,f_B)$\\
		($(A,B)\in CD_X^P(f_A,f_B)$).
	\end{defn}
	
	Let $(X,\rho)$ be a metric space, $A$ and $B$ be subsets of $X$, $C$ be an arbitrary set, and $f:C\to \mathbb{R},\ g:C\to \mathbb{R}$. Wi will use the notation\\
	$\displaystyle S_{f,g}^C(A,B)= \dist(A,B)+\inf_{c\in C}f(c)+\inf_{c\in C}g(c)$
	
	With the following definition, we will introduce a contraction condition that the basic theorem of this article is related to.
	
	\begin{defn}\label{contraction}
		Let $(X,\rho)$ be a metric space, $A$ and $B$ be subsets of $X$, $C$ be an arbitrary set, $P\subseteq A\times B \times C^2$, and
		$T_A:A\times C \to A,\ H_A:A\times C \to C,\ T_B:B\times C \to B,\ H_B:B\times C \to C,\ f_A:
		C\to \mathbb{R},\ f_B:C\to \mathbb{R}$ be such that
		\begin{enumerate}[label=C-\arabic*:]
			\item\label{C-1} If $(x_0,y_0,u_0,v_0)\in P$, $\{x_n\}_{n=0}^{\infty}$ and $\{u_n\}_{n=0}^{\infty}$ are the iterated sequences of $(T_A,H_A)$ with an initial guess $(x_0,u_0)$, and $\{y_n\}_{n=0}^{\infty}$ and $\{v_n\}_{n=0}^{\infty}$ are the iterated sequences of $(T_B,H_B)$ with an initial guess $(y_0,v_0)$, then the inclusion $(x_n,y_m,u_n,v_m)\in P$ holds for every $n,m\in \mathbb{N}$. 
			\item there holds
			\begin{equation}\label{d1e2}
				\inf_{c\in C}f_A(c)\in \mathbb{R}\ and\ \inf_{c\in C}f_B(c)\in \mathbb{R}
			\end{equation}
			\item there is $\lambda\in[0,1)$ so that for every $(x,y,u,v)\in P$ the inequality
			\begin{equation}\label{d1e1}
				\begin{array}{rcl}
					W_2&=&\rho(T_A(x,u),T_B(y,v))+f_A(H_A(x,u))+f_B(H_B(y,v))\\
					&\leq&\lambda(\rho(x,y)+f_A(u)+f_B(v))+(1-\lambda)S_{f_A,f_B}^C(A,B)
				\end{array}
			\end{equation}
			is true.
		\end{enumerate}
		Then we say that the ordered set $(T_A,T_B,H_A,H_B,f_A,f_B)$ is a contraction map set with an external factor about the sets $(A,B,C) $ in $P$. Forward we will notate this with $(T_A,T_B,H_A,H_B,f_A,f_B)\in CEF_C^P(A,B)$.
	\end{defn}
	
	The next definition is technical.
	
	\begin{defn}\label{def6}
		Let $A$, $B$, and $C$ be sets, $a\in A$, $P\subseteq A\times B \times C^2$, and $ f:
		C\to \mathbb{R}$, let the sequence $\{c_n\}_{n=0}^{\infty}\subset C$ be such that
		\begin{itemize}
			\item $\displaystyle\lim_{n\to \infty}f(c_n)=\inf_{c\in C}f(c)$
			\item there exist $y\in B$ and $v\in C$ so that $(a,y,c_n,v)\in P$ for every $n\in \mathbb{N}$
		\end{itemize}
		Then we say that $\{c_n\}_{n=0}^{\infty}$ is an infimum sequence of $a$ on $f$ and $(y,v)$ in $P$ ($\{c_n\}_{n=0}^{\infty}\in IS_f^P(a,y,v)$).
	\end{defn} 
	
	Now we will formulate the theorem on which the main focus of this article falls.
	
	\begin{thm}\label{t1}
		Let $(X,\rho)$ be a metric space, $A$ and $B$ be subsets of $X$, $C$ be an arbitrary set, $P\subseteq A\times B\times C^2$, $(T_A,T_B,H_A,H_B,f_A,f_B)\in CEF_C^P(A,B)$, $(A,B)\in CD_X^P(f_A,f_B)$. Then:
		\begin{enumerate}
			\item\label{t1-cons1} If $(x_0,y_0,u_0,v_0)\in P$, $\{x_n\}_{n=0}^\infty$ and $\{u_n\}_{n=0}^\infty$ are the iterated sequences of $(T_A,H_A)$ with an initial guess $(x_0,u_0)$, and $\{y_n\}_{n=0}^\infty$ and $\{v_n\}_{n=0}^\infty$ are the iterated sequences of $(T_B,H_B)$ with an initial guess $(y_0,v_0)$, then there exists $\displaystyle\lim_{n\to\infty}x_n=\alpha\in A$ and
			there hold
			$$
			\begin{array}{l}
				\displaystyle\lim_{n\to \infty}\rho(y_n,\alpha)= \dist(A,B)\\
				\displaystyle\lim_{n\to \infty} f_A(u_n)= \inf_{c\in C} f_A(c)\\
				\displaystyle\lim_{n\to \infty}f_B(v_n)= \inf_{c\in C} f_B(c)
			\end{array}
			$$
			\item\label{t1-cons2} If the inclusions $(x_0,y_0,u_0.v_0)\in P$ and $(z_0,y_0,t_0,v_0)\in P$ hold, $\{x_n\}_{n=0}^\infty$ and $\{u_n\}_{n=0}^\infty$ are the iterated sequences of $(T_A,H_A)$ with an initial guess $(x_0,u_0)$, and $\{z_n\}_{n=0}^\infty$ and $\{t_n\}_{n=0}^\infty$ are the iterated sequences of $(T_A,H_A)$ with an initial guess $(z_0,t_0)$, then $\displaystyle\lim_{n\to \infty}x_n=\lim_{n\to \infty}z_n\in A$.
			\item\label{t1-cons3} If $(x_0,y_0,u_0,v_0)\in P$, $\{x_n\}_{n=0}^\infty$ and $\{u_n\}_{n=0}^\infty$ are the iterated sequences of $(T_A,H_A)$ with an initial guess $(x_0,u_0)$, $\alpha=\displaystyle\lim_{n\to \infty}x_{n}$ and\\
			$\{c_n\}_{n=0}^\infty\in IS_{f_A}^P(\alpha,y_0,v_0)$. Then $\displaystyle\lim_{n\to \infty}T_A(\alpha,c_n)=\alpha$.
			\item\label{t1-cons4} If $(x_0,y_0,u_0,v_0)\in P$, $\{x_n\}_{n=0}^\infty$ and $\{u_n\}_{n=0}^\infty$ are the iterated sequences of $(T_A,H_A)$ with an initial guess $(x_0,u_0)$, and $\displaystyle\lim_{n\to\infty}x_n=\alpha$. Then there is no $\beta\in A$ such that $\beta\neq \alpha$ and $\displaystyle\lim_{n\to \infty}T_A(\beta,c_n)=\beta$, where $\{c_n\}_{n=0}^\infty\in IS_{f_A}^P(\beta,y_0,v_0)$.
		\end{enumerate}
	\end{thm}
	
	\section{\normalsize Auxiliary Results}
	
	\subsection{\normalsize Some auxiliary results about convergent sequences}
	
	We will start with some results about convergent sequences in metric spaces.
	
	\begin{lem}\label{l0}
		Let $\{x_n\}_{n=0}^\infty\subset \mathbb{R}$ and $\{y_n\}_{n=0}^\infty\subset \mathbb{R}$. Let for every $\varepsilon>0$, there exists $N\in \mathbb{N}$, so that if $n\geq N$ then $\displaystyle x_n+y_n\leq x'+y'+\varepsilon$ holds. Where $\displaystyle x'\leq\inf_{n\in \mathbb{N}}x_n$ and $\displaystyle y'\leq \inf_{n\in \mathbb{N}}y_n$. Then $\displaystyle \lim_{n\to\infty}x_n=x'$ and $\displaystyle\lim_{n\to\infty}y_n=y'$.
	\end{lem}
	\begin{proof}
		From $\displaystyle x_n+y_n\leq x'+y'+\varepsilon$ and $\displaystyle y'\leq\inf_{n\in \mathbb{N}}y_n\leq y_n$ we get $\displaystyle x'\leq\inf_{n\in \mathbb{N}}x_n\leq x_n\leq x'+\varepsilon$, i.e., for every $\varepsilon>0$ there exists $N\in \mathbb{N}$ such that if $n\geq N$ then  $\displaystyle x'\leq x_n\leq x'+\varepsilon$ is true. Thus, $\displaystyle \lim_{n\to\infty}x_n=x'$. Similarly, we can prove $\displaystyle\lim_{n\to\infty}y_n=y'$.
	\end{proof}
	Similarly to lemma \ref{l0}, we can prove the next corollary.
	\begin{cor}\label{c1}
		If $\{x_n\}_{n=0}^\infty\subset \mathbb{R}$, $\{y_n\}_{n=0}^\infty\subset \mathbb{R}$, $\{z_n\}_{n=0}^\infty\subset \mathbb{R}$ and for every $\varepsilon>0$ there exists $N\in \mathbb{N}$ so that for any $n\geq N$ there holds $\displaystyle x_n+y_n+z_n\leq x'+y'+z'+\varepsilon$. Where $\displaystyle x'\leq\inf_{n\in \mathbb{N}}x_n$, $\displaystyle y'\leq\inf_{n\in \mathbb{N}}y_n$, and $\displaystyle z'\leq\inf_{n\in \mathbb{N}}z_n$. Then $\displaystyle \lim_{n\to\infty}x_n=x'$, $\displaystyle \lim_{n\to\infty}y_n=y'$, and $\displaystyle \lim_{n\to\infty}z_n=z'$.
	\end{cor}
	\begin{lem}\label{l0.0}
		Let $f:\mathbb{N}^2\to \mathbb{R}$, $g:\mathbb{N}^2\to \mathbb{R}$, and 
		$$\displaystyle\lim_{k\to \infty}\sup_{n\geq k,\ m\geq k} (f(n,m)+g(n,m))=\inf_{n,m\in \mathbb{N}}f(n,m)+\inf_{n,m\in \mathbb{N}}g(n,m)>-\infty.$$
		Then  
		\begin{itemize}
			\item $\displaystyle\lim_{k\to \infty}\sup_{n\geq k,\ m\geq k}f(n,m)=\inf_{n,m\in \mathbb{N}}f(n,m)$
			\item $\displaystyle\lim_{k\to \infty}\sup_{n\geq k,\ m\geq k}g(n,m)=\inf_{n,m\in \mathbb{N}}g(n,m).$
		\end{itemize}
		
	\end{lem}
	\begin{proof}
		Let us suppose that $\displaystyle\lim_{k\to \infty}\sup_{n\geq k,\ m\geq k}f(n,m)\neq\inf_{n,m\in \mathbb{N}}f(n,m)$. Thus there is $\varepsilon>0$ such that for every $k\in \mathbb{N}$ there exist $n,m\geq k$ so that $\displaystyle f(n,m)>\inf_{n,m\in \mathbb{N}}f(n,m)+\varepsilon$. Therefore for any $k\in \mathbb{N}$ there holds $\displaystyle\sup_{n\geq k,\ m\geq k}(f(n,m)+g(n,m))>\inf_{n,m\in \mathbb{N}}f(n,m)+\inf_{n,m\in \mathbb{N}}g(n,m)+\varepsilon$. From the last inequality and
		$\displaystyle\inf_{n,m\in \mathbb{N}}f(n,m)+\inf_{n,m\in \mathbb{N}}g(n,m)\in \mathbb{R} $ we get
		$$\displaystyle\lim_{k\to \infty}\sup_{n\geq k,\ m\geq k}(f(n,m)+g(n,m))\neq\inf_{n,m\in \mathbb{N}}f(n,m)+\inf_{n,m\in \mathbb{N}}g(n,m).$$
		This is a contradiction.
		
		Similarly, we can prove that if $\displaystyle\lim_{k\to \infty}\sup_{n\geq k,\ m\geq k}g(n,m)\neq\inf_{n,m\in \mathbb{N}}g(n,m)$, a contradiction also arises.
	\end{proof}
	\begin{cor}\label{c1.1}
		Let $\{x_n\}_{n=0}^\infty\subset \mathbb{R}$, $\{y_n\}_{n=0}^\infty\subset \mathbb{R}$, $f:\mathbb{N}^2\to \mathbb{R}$, and 
		$$\displaystyle\lim_{k\to \infty}\sup_{n\geq k,\ m\geq k} (f(n,m)+x_n+y_m)=\inf_{n,m\in \mathbb{N}}f(n,m)+\inf_{n\in \mathbb{N}}x_n+\inf_{n\in \mathbb{N}}y_n\in \mathbb{R}.$$ Then
		$$\lim_{k\to \infty}\sup_{n\geq k,\ m\geq k}f(n,m)=\inf_{n,m\in \mathbb{N}}f(n,m).$$
	\end{cor}
	\begin{proof}
		We can denote $g(n,m)=x_n+y_m$ and $\displaystyle \inf_{n\in \mathbb{N}}x_n+\inf_{n\in \mathbb{N}}y_n=\inf_{n,m\in \mathbb{N}}g(n,m)$. Then this corollary follows from lemma \ref{l0.0}.
	\end{proof}
	\subsection{\normalsize The property Convergent By distance}
	
	\begin{lem}\label{l0.1}
		Let $(X,\rho)$ be a metric space, $A$ and $B$ be subsets of $X$, $C$ be an arbitrary set, $f_A:C\to \mathbb{R}$ and $f_B:C\to \mathbb{R}$. Let $(A,B)\in CD_X^P(f_A,f_B)$, the sequences $\{x_n\}_{n=0}^\infty\subseteq A$, $\{z_n\}_{n=0}^\infty\subseteq A$, and $\{y_n\}_{n=0}^\infty\subseteq B$ be such that
		\begin{itemize}
			\item \begin{equation}\label{l0.1e0} 
				\left|
				\begin{array}{l}
					\displaystyle\lim_{k\to \infty}\sup_{n\geq k,\ m\geq k}\rho(x_n,y_m)= \dist(A,B)\\
					\displaystyle\lim_{k\to \infty}\sup_{n\geq k,\ m\geq k}\rho(z_n,y_m)= \dist(A,B)
				\end{array}
				\right.
			\end{equation}
			\item there exist sequences $\{u_n\}_{n=1}^{\infty}\subseteq C$, $\{v_n\}_{n=1}^{\infty}\subseteq C$, and $\{t_n\}_{n=1}^{\infty}\subseteq C$ so that 
			\begin{itemize}
				\item[${\scriptscriptstyle\blacklozenge}$] $(x_n,y_m,u_n,v_m)\in P$ for every $n,m\in \mathbb{N}$
				\item[${\scriptscriptstyle\blacklozenge}$]$(z_n,y_m,t_n,v_m)\in P$ for every $n,m\in \mathbb{N}$
				\item[${\scriptscriptstyle\blacklozenge}$] $\displaystyle \lim_{n\to \infty}f_A(u_n)=\inf_{c\in C}f_A(c) $, $\displaystyle \lim_{n\to \infty}f_B(v_n)=\inf_{c\in C}f_B(c) $ and\\
				$\displaystyle \lim_{n\to \infty}f_A(t_n)=\inf_{c\in C}f_A(c) $
			\end{itemize}
		\end{itemize}
		then there exists $\alpha\in A$ such that $\alpha=\displaystyle\lim_{n\to \infty}x_n=\lim_{n\to \infty}z_n$.
	\end{lem}
	\begin{proof}
		Let the sequences $\{a_n\}_{n=1}^\infty$, $\{b_n\}_{n=1}^\infty$, $\{{a_C}_n\}_{n=1}^\infty$, and $\{{b_C}_n\}_{n=1}^\infty$ be defined as
		$a_{2q-1}=x_q$, $a_{2q}=z_q$, ${a_C}_{2q-1}=u_q$, ${a_C}_{2q}=t_q$ for every $q\in \mathbb{N}$, $b_{n}=y_{\lceil\frac{n}{2}\rceil}$ and ${b_C}_{n}=v_{\lceil\frac{n}{2}\rceil}$. We can see
		\begin{equation}\label{l0.1e0.1} 
			\displaystyle \lim_{n\to \infty}f_A({a_C}_n)=\inf_{c\in C}f_A(c)\ and\ \displaystyle \lim_{n\to \infty}f_B({b_C}_n)=\inf_{c\in C}f_B(c).
		\end{equation}
		For $\{a_n\}_{n=1}^\infty\subseteq A$ and $\{b_n\}_{n=1}^\infty\subseteq B$, we can also see that for every $k\in \mathbb{N}$, the equality
		\begin{equation}\label{l0.1e1} 
			\begin{array}{rcl}
				W_3&=&\ds\sup_{m\geq k,\ n\geq k}\rho(b_n,a_m)\\[20pt]
				&=&\ds\max\left\{\sup_{n\geq k,\ 2q\geq k}\rho(b_{n},a_{2q}),\sup_{n\geq k,\ 2q-1\geq k}\rho(b_{n},a_{2q-1})\right\}
			\end{array}
		\end{equation}
		holds. By the definitions of $\{b_n\}_{n=1}^\infty$ and $\{{a_C}_n\}_{n=1}^\infty$, we get for every $k\in \mathbb{N}$ there holds
		$$
		\begin{array}{rcl}
			\dist(A,B)&\leq&\ds\sup_{n\geq k,\ 2q\geq k}\rho(b_{n},a_{2q})=\sup_{n\geq k,\ 2q \geq k}\rho(y_{\lceil\frac{n}{2}\rceil},z_q)\\[20pt]
			&\leq&\ds\sup_{\lceil\frac{n}{2}\rceil\geq \lceil\frac{k}{2}\rceil,\ q\geq \lceil\frac{k}{2}\rceil}\rho(y_{\lceil\frac{n}{2}\rceil},z_q)=\sup_{i\geq \lceil\frac{k}{2}\rceil,\ m\geq \lceil\frac{k}{2}\rceil}\rho(z_i,y_m).
		\end{array}
		$$
		By the last inequality and (\ref{l0.1e0}), it follows
		\begin{equation}\label{l0.1e2} 
			\lim_{k\to \infty}\sup_{n\geq k,\ 2q\geq k}\rho(b_{n},a_{2q})= \dist(A,B).
		\end{equation} 
		Using the definitions of $\{b_n\}_{n=1}^\infty$ and $\{{a_C}_n\}_{n=1}^\infty$, we obtain that for every $k\in \mathbb{N}$, the chain of inequalities
		$$ 
		\begin{array}{rcl}
			\dist(A,B)&\leq&\displaystyle\sup_{n\geq k,\ 2q-1\geq k}\rho(b_{n},a_{2q-1})\leq\sup_{\lceil\frac{n}{2}\rceil\geq \lceil\frac{k}{2}\rceil,\ 2q>k}\rho(y_{\lceil\frac{n}{2}\rceil},x_{q})\\[20pt]
			&\leq&\ds\sup_{\lceil\frac{n}{2}\rceil\geq \lceil\frac{k}{2}\rceil,\ q\geq \lceil\frac{k}{2}\rceil}\rho(y_{\lceil\frac{n}{2}\rceil},x_{q})=\sup_{i\geq \lceil\frac{k}{2}\rceil,\ m\geq \lceil\frac{k}{2}\rceil}\rho(x_i,y_m)
		\end{array}
		$$
		holds. From the above one and (\ref{l0.1e0}), we can observe that
		$$
		\lim_{k\to \infty}\sup_{n\geq k,\ 2q-1\geq k}\rho(b_{n},a_{2q-1})= \dist(A,B).
		$$
		By the last equality, (\ref{l0.1e2}) and (\ref{l0.1e1}), we get
		\begin{equation}\label{l0.1e3} 
			\lim_{k\to \infty}\sup_{m\geq k,\ n\geq k}\rho(b_n,a_m)= \dist(A,B).
		\end{equation} 
		We can see that every $q,m\in \mathbb{N}$ there holds 
		$$
		\begin{array}{l}
			(a_{2q-1},b_m,{a_C}_{2q-1},{b_C}_m)=\left(x_q,y_{\lceil\frac{m}{2}\rceil},u_q,v_{\lceil\frac{m}{2}\rceil}\right)\in P\\
			(a_{2q},b_m,{a_C}_{2q},{b_C}_m)=\left(z_q,y_{\lceil\frac{m}{2}\rceil},t_q,v_{\lceil\frac{m}{2}\rceil}\right)\in P,
		\end{array}
		$$
		i.e., for every $n,m\in \mathbb{N}$ there holds $(a_{n},b_m,{a_C}_{n},{b_C}_m)\in P$.
		Using the last inclusion, (\ref{l0.1e3}), (\ref{l0.1e0.1}), and $(A,B)\in CD_X^P(f_A,f_B)$, it follows that
		$$
		\lim_{q\to \infty}x_{q}=\lim_{q\to \infty}a_{2q-1}=\lim_{n\to \infty}a_n=\lim_{q\to \infty}a_{2q}=\lim_{q\to \infty}z_{q}\ and\ \lim_{n\to \infty}a_n\in A.
		$$
	\end{proof} 
	\subsection{\normalsize A contraction map set with an external factor about sets}
	
	\begin{lem}\label{l1}
		Let $(X,\rho)$ be a metric space, $A, B\subseteq X$, $C$ be an arbitrary set, $P\subseteq A\times B\times C^2$, and $(T_A,T_B,H_A,H_B,f_A,f_B)\in CEF_C^P(A,B)$. Let $(x_0,y_0,u_0,v_0)\in P$, $\{x_n\}_{n=0}^\infty$ and $\{u_n\}_{n=0}^\infty$ be the iterated sequences of\\
		$(T_A,H_A)$ with an initial guess $(x_0,u_0)$, and  $\{y_n\}_{n=0}^\infty$ and $\{v_n\}_{n=0}^\infty$ be the iterated sequences of $(T_B,H_B)$ with an initial guess $(y_0,v_0)$. Then
		\begin{itemize}
			\item $\{x_n\}_{n=0}^\infty$ and $\{y_n\}_{n=0}^\infty$ are bounded
			\item $\displaystyle \sup_{n\in \mathbb{N}}f_A(u_n)<\infty$ and $\displaystyle\sup_{n\in \mathbb{N}}f_B(v_n)<\infty.$
		\end{itemize}
	\end{lem}
	\begin{proof}
		From $(T_A,T_B,H_A,H_B,f_A,f_B)\in CEF_C^P(A,B)$, definition \ref{def4} and definition \ref{contraction}:\ref{C-1} we get, for every $n\in \mathbb{N}$ there holds
		\begin{equation}\label{l1e1}
			\begin{array}{rcl}
				W_4&=&\rho(x_{n+1},y_2)+f_A(u_{n+1})+f_B(v_2)\\
				&\leq&\lambda(\rho(x_n,y_1)+f_A(u_n)+f_B(v_1))+(1-\lambda)S_{f_A,f_B}^C(A,B)
			\end{array}
		\end{equation}
		Using the triangle inequality, we have $\rho(x_{n+1},y_2)\geq \rho(x_{n+1},y_1)-\rho(y_1,y_2)$. By the last inequality and (\ref{l1e1}), it follows
		\begin{equation}\label{l1e2}
			\begin{array}{rcl}
				W_5&=&\rho(x_{n+1},y_1)+f_A(u_{n+1})\\
				&\leq&\lambda(\rho(x_n,y_1)+f_A(u_n))+Q+(1-\lambda)S_{f_A,f_B}^C(A,B)
			\end{array}
		\end{equation}
		where $Q=\rho(y_1,y_2)+\lambda f_B(v_1)-f_B(v_2)$. Applying (\ref{l1e2}) $n$ times consecutively, we get
		$$
		\begin{array}{rl}
			W_6=&\rho(x_n,y_1)+f_A(u_n)\\
			\leq&\lambda(\rho(x_{n-1},y_1)+f_A(u_{n-1}))+Q+(1-\lambda)S_{f_A,f_B}^C(A,B)\\
			\leq&\lambda^2(\rho(x_{n-2},y_1)+f_A(u_{n-2}))+Q+\lambda Q+(1-\lambda^2)S_{f_A,f_B}^C(A,B)\\
			\dots\\
			\leq&\displaystyle\lambda^{n-3}(\rho(x_{3},y_1)+f_A(u_{3}))+Q\frac{1-\lambda^{n-3}}{1-\lambda}+(1-\lambda^{n-3})S_{f_A,f_B}^C(A,B)\\[12pt]
			\leq&\displaystyle\lambda^{n-2}(\rho(x_{2},y_1)+f_A(u_{2}))+Q\frac{1-\lambda^{n-2}}{1-\lambda}+(1-\lambda^{n-2})S_{f_A,f_B}^C(A,B)\\[12pt]
			\leq&\displaystyle\lambda^{n-1}(\rho(x_{1},y_1)+f_A(u_{1}))+Q\frac{1-\lambda^{n-1}}{1-\lambda}+(1-\lambda^{n-1})S_{f_A,f_B}^C(A,B)\\[12pt]
		\end{array}
		$$
		and hence, for every $n\in \mathbb{N}$, the inequality
		$$\rho(x_n,y_1)+f_A(u_n)\leq \rho(x_{1},y_1)+f_A(u_{1})+Q\frac{1}{1-\lambda}+S_{f_A,f_B}^C(A,B)$$
		holds. Thus
		$$\sup_{n\in \mathbb{N}}\{\rho(x_n,y_1)+f_A(u_n)\}< \infty.$$ 
		From the last inequality and (\ref{d1e2}), we can observe that $\displaystyle\sup_{n\in \mathbb{N}}f_A(u_n)< \infty$ and $\displaystyle\sup_{n\in \mathbb{N}}\{\rho(x_n,y_1)\}<\infty$, i.e., $\{x_n\}_{n=0}^\infty$ is a bounded sequence. By similar arguments, it follows that $\displaystyle\sup_{n\in \mathbb{N}}f_B(v_n)< \infty$ and that $\{y_n\}_{n=0}^\infty$ is a bounded sequence.
	\end{proof}
	\begin{lem}\label{l2}
		Let $(X,\rho)$ be a metric space, $A$ and $B$ be subsets of $X$, $C$ be an arbitrary set, $P\subseteq A\times B\times C^2$, $(A,B)\in CD_X^P(f_A,f_B)$, and\\
		$(T_A,T_B,H_A,H_B,f_A,f_B)\in CEF_C^P(A,B)$. Let $(x_0,y_0,u_0,v_0)\in P$, $\{x_n\}_{n=0}^\infty$ and $\{u_n\}_{n=0}^\infty$ be the iterated sequences of $(T_A,H_A)$ with an initial guess $(x_0,u_0)$,  $\{y_n\}_{n=0}^\infty$ and $\{v_n\}_{n=0}^\infty$ be the iterated sequences of $(T_B,H_B)$ with an initial guess $(y_0,v_0)$ then:
		\begin{itemize}
			\item $\{x_n\}_{n=0}^\infty$ converges in $A$
			\item $\displaystyle\lim_{k\to \infty} \sup_{n\geq k,m\geq k}\rho(x_n,y_{m})=\lim_{n\to \infty} \rho(x_n,y_{n})=\lim_{n\to \infty} \rho(\alpha,y_{n})= \dist(A,B)$ where $\displaystyle\lim_{n\to \infty}x_{n}=\alpha$
			\item $\displaystyle\lim_{n\to \infty} f_A(u_n)= \inf_{c\in C} f_A(c)$ and $\displaystyle \lim_{n\to \infty}f_B(v_n)= \inf_{c\in C} f_B(c).$
		\end{itemize}
	\end{lem}
	\begin{proof}
		From $\{x_n\}_{n=0}^\infty\subseteq A$ and $\{y_n\}_{n=0}^\infty\subseteq B$, we get 
		\begin{equation}\label{l2e0}
			\displaystyle\dist(A,B)\leq\inf_{n,m\in\mathbb{N}}\rho(x_{m},y_{n}).
		\end{equation}
		
		By $(T_A,T_B,H_A,H_B,f_A,f_B)\in CEF_C^P(A,B)$, definition \ref{def4} and definition \ref{contraction}:\ref{C-1} it follows that for every $n\in \mathbb{N}$ and $m\in \mathbb{N}$ so that $n\leq m$ the inequality
		$$
		\begin{array}{rl}
			W_7=&\rho(x_{m},y_{n})+f_A(u_{m})+f_B(v_{n})\\
			\leq&\lambda(\rho(x_{m-1},y_{n-1})+f_A(u_{m-1})+f_B(v_{n-1}))+(1-\lambda)S_{f_A,f_B}^C(A,B)
		\end{array}
		$$
		holds. Applying the above one $n-1$ times consecutively, we obtain
		\begin{equation}\label{l2e1}
			\begin{array}{rcl}
				U(m,n)&\leq&\lambda U(m-1,n-1)+(1-\lambda)S_{f_A,f_B}^C(A,B)\\
				&\leq&\lambda^2 U(m-2,n-2)+(1-\lambda^2)S_{f_A,f_B}^C(A,B)\\
				&\leq&\lambda^3 U(m-3,n-3)+(1-\lambda^3)S_{f_A,f_B}^C(A,B)\\
				&\dots\\
				&\leq&\lambda^{n-3} U(m-n+3,3)+(1-\lambda^{n-3})S_{f_A,f_B}^C(A,B)\\
				&\leq&\lambda^{n-2} U(m-n+2,2)+(1-\lambda^{n-2})S_{f_A,f_B}^C(A,B)\\
				&\leq&\lambda^{n-1} U(m-n+1,1)+(1-\lambda^{n-1})S_{f_A,f_B}^C(A,B)\\
			\end{array}
		\end{equation}
		where $U(m,n)=\rho(x_{m},y_{n})+f_A(u_{m})+f_B(v_{n})$. Let 
		$$\displaystyle M=\max\left\{\sup_{k\in \mathbb{N}}U(k,1),\sup_{k\in \mathbb{N}}U(1,k)\right\}.$$ From lemma \ref{l1}, we have $M<\infty$.
		Using $M<\infty$ and (\ref{l2e1}), we can observe that for every $n\in \mathbb{N}$ and $m\in \mathbb{N}$ such that $n\leq m$, the inequality 
		\begin{equation}\label{l2e1.1}
			\rho(x_{m},y_{n})+f_A(u_{m})+f_B(v_{n})\leq\lambda^{n-1}M+(1-\lambda^{n-1})S_{f_A,f_B}^C(A,B)
		\end{equation}
		is true.
		By similar arguments, we get that for every $n\in \mathbb{N}$ and $m\in \mathbb{N}$ so that $n\geq m$, there holds
		\begin{equation}\label{l2e1.2}
			\rho(x_{m},y_{n})+f_A(u_{m})+f_B(v_{n})\leq\lambda^{m-1}M+(1-\lambda^{m-1})S_{f_A,f_B}^C(A,B).
		\end{equation}
		From (\ref{l2e1.1}), it follows that for every $\varepsilon>0$, there exists $N\in \mathbb{N}$ such that if $N\leq n$, then $\rho(x_{n},y_{n})+f_A(u_{n})+f_B(v_{n})\leq S_{f_A,f_B}^C(A,B)+\varepsilon$ is true. By the last inequality, (\ref{l2e0}) and corollary \ref{c1}, we obtain
		\begin{equation}\label{l2e2}
			\lim_{n\to \infty} \rho(x_{n},y_{n})= \dist(A,B)
		\end{equation}
		\begin{equation}\label{l2e3}
			\displaystyle\lim_{n\to \infty} f_A(u_n)= \inf_{c\in C} f_A(c),\ \displaystyle \lim_{n\to \infty}f_B(v_n)= \inf_{c\in C} f_B(c).
		\end{equation}
		
		Using $\{x_n\}_{n=0}^\infty\subseteq A$, $\{y_n\}_{n=0}^\infty\subseteq B$, and (\ref{l2e2}), we get 
		\begin{equation}\label{l2e4.1}
			\displaystyle\inf_{n,m\in\mathbb{N}}\rho(x_{m},y_{n})= \dist(A,B).
		\end{equation}
		
		From (\ref{l2e1.1}), (\ref{l2e1.2}), and (\ref{l2e4.1}), we can observe that for every $n\in \mathbb{N}$ and $m\in \mathbb{N}$ there holds 
		$$
		\begin{array}{rcl}
			S_{f_A,f_B}^C(A,B)&\leq&\ds\rho(x_{m},y_{n})+f_A(u_{m})+f_B(v_{n})\\[10pt]
			&\leq&\ds\lambda^{\min\{m,n\}-1}M+(1-\lambda^{\min\{m,n\}-1})S_{f_A,f_B}^C(A,B).
		\end{array}
		$$
		Thus we get
		\begin{equation}\label{l2e4}
			\displaystyle\lim_{k\to \infty}\sup_{n\geq k,\ m\geq k}(\rho(x_m,y_n)+f_A(u_{m})+f_B(v_{n}))=S_{f_A,f_B}^C(A,B).
		\end{equation}
		
		By $(T_A,T_B,H_A,H_B,f_A,f_B)\in CEF_C^P(A,B)$, we get
		$\displaystyle \inf_{n\in \mathbb{N}}f_A(u_n)>-\infty$ and $\displaystyle\inf_{n\in \mathbb{N}}f_B(v_n)>-\infty$. 
		By the last two inequalities, (\ref{l2e4}), (\ref{l2e4.1}), and Corollary \ref{c1.1}, it follows
		\begin{equation}\label{l2e5}
			\displaystyle\lim_{k\to \infty}\sup_{n\geq k,\ m\geq k}\rho(x_m,y_n)= \dist(A,B).
		\end{equation}
		From definition \ref{contraction}:\ref{C-1} we get that for every $n,m\in \mathbb{N}$ there holds\\
		$(x_n,y_m,u_n,v_m)\in P$. Using the last inclusion, (\ref{l2e3}), (\ref{l2e5}) and $(A,B)\in CD_X^P(f_A,f_B)$, we can observe that the sequence $\{x_n\}_{n=1}^\infty$ converges in $A$.
		
		Let $\displaystyle\lim_{n\to \infty}x_{n}=\alpha$.Then by the continuity of the metric function and (\ref{l2e2}) it follows $\displaystyle\lim_{n\to \infty} \rho(\alpha,y_{n})= \dist(A,B)$.
	\end{proof}
	\begin{lem}\label{l2.1}
		Let $(X,\rho)$ be a metric space, $A$ and $B$ be subsets of $X$,\\
		$C$ be an arbitrary set, $P\subseteq A\times B\times C^2$, $(A,B)\in CD_X^P(f_A,f_B)$, and $(T_A,T_B,H_A,H_B,f_A,f_B)\in CEF_C^P(A,B)$. Let $(x_0,y_0,u_0.v_0)\in P$ and\\
		$(z_0,y_0,t_0,v_0)\in P$ Then $\displaystyle\lim_{n\to \infty}x_n=\lim_{n\to \infty}z_n\in A$ where $\{x_n\}_{n=0}^\infty$, $\{u_n\}_{n=0}^\infty$ are the iterated sequences of $(T_A,H_A)$ with an initial guess $(x_0,u_0)$ and $\{z_n\}_{n=0}^\infty$, $\{t_n\}_{n=0}^\infty$ are the iterated sequences of $(T_A,H_A)$ with an initial guess $(z_0,t_0)$.
	\end{lem}
	\begin{proof}
		Let $\{y_n\}_{n=0}^\infty\subset B$ and $\{v_n\}_{n=0}^\infty\subset C$ be the iterated sequences of $(T_B,H_B)$ with an initial guess $(y_0,v_0)$.
		From lemma \ref{l2} we get that there exist $\displaystyle\lim_{n\to \infty}x_n=\alpha$, $\displaystyle\lim_{n\to \infty}z_n=\beta$, and
		\begin{equation}\label{l2.1e1}
			\begin{array}{rcl}
				\ds\lim_{n\to \infty}f_A(u_n)&=&\ds\inf_{c\in C} f_A(c)\\[10pt]
				\ds\lim_{n\to \infty}f_B(v_n)&=&\ds\inf_{c\in C} f_B(c)\\[10pt]
				\ds\lim_{n\to \infty}f_A(t_n)&=&\ds\inf_{c\in C} f_A(c).
			\end{array}
		\end{equation}
		By lemma \ref{l2}, it follows
		$$
		\left|
		\begin{array}{l}
			\displaystyle \lim_{k\to \infty} \sup_{n\geq k,m\geq k}\rho(x_n,y_{m})= \dist(A,B)\\
			\displaystyle \lim_{k\to \infty} \sup_{n\geq k,m\geq k}\rho(z_n,y_{m})= \dist(A,B).
		\end{array}
		\right.
		$$
		Using the last system, lemma \ref{l0.1}, definition \ref{contraction}:\ref{C-1}, (\ref{l2.1e1}) and $(A,B)\in CD_X^P(f_A,f_B)$, we obtain $\alpha=\beta$
	\end{proof}
	
	\subsection{\normalsize A weakly fixed point of a map about the sequence}
	
	With the following definition, we extend the concept of a fixed point, introducing the concept of a weak fixation of a point about a sequence.
	
	\begin{lem}\label{l3}
		Let $(X,\rho)$ be a metric space, $A$ and $B$ be subsets of $X$, $C$ be an arbitrary set, $P\subseteq A\times B\times C^2$, $(A,B)\in CD_X^P(f_A,f_B)$, and\\
		$(T_A,T_B,H_A,H_B,f_A,f_B)\in CEF_C^P(A,B)$. Let $(x_0,y_0,u_0,v_0)\in P$, $\{x_n\}_{n=0}^\infty$, $\{u_n\}_{n=0}^\infty$ be the iterated sequences of $(T_A,H_A)$ with an initial guess $(x_0,u_0)$, $\alpha=\displaystyle\lim_{n\to \infty}x_{n}$, and $\{c_n\}_{n=0}^\infty\in IS_{f_A}^P(\alpha,y_0,v_0)$. Then $\displaystyle\lim_{n\to \infty}T_A(\alpha,c_n)=\alpha$.
	\end{lem}
	\begin{proof}
		Let $\{y_n\}_{n=0}^\infty$ and $\{v_n\}_{n=0}^\infty$ be the iterated sequences of $(T_B,H_B)$ with an initial guess $(y_0,v_0)$.
		For every $n,m\in \mathbb{N}$, there holds $y_{n+1}\in B$ and $T_A(\alpha,c_m)\in A$, therefore
		\begin{equation}\label{l3e1.3} 
			\rho(y_{n+1},T_A(\alpha,c_m))\geq  \dist(A,B).
		\end{equation} 
		From definition \ref{contraction}:\ref{C-1},definition \ref{def4}, and $\{c_n\}_{n=0}^\infty\in IS_{f_A}^P(\alpha,y_0,v_0)$, we can observe that for any $n,m\in \mathbb{N}$ the inclusion
		\begin{equation}\label{l3e-1}
			(\alpha,y_m,c_n,v_m)\in P
		\end{equation}  
		is true. By (\ref{l3e1.3}), (\ref{l3e-1}), and $(T_A,T_B,H_A,H_B,f_A,f_B)\in CEF_C^P(A,B)$, we can see that for every $n,m\in \mathbb{N}$, there holds
		\begin{equation}\label{l3e0}
			\begin{array}{rl}
				S\leq&\rho(y_{n+1},T_A(\alpha,c_m))+I\\
				\leq&f_A(H_A(\alpha,c_m))+\rho(y_{n+1},T_A(\alpha,c_m))+f_B(v_{n+1})\\
				\leq&\lambda(f_A(c_m)+\rho(y_{n},\alpha)+f_B(v_{n}))+(1-\lambda)S_{f_A,f_B}^C(A,B),
			\end{array}
		\end{equation}
		where $\displaystyle I=\inf_{c\in C}f_A(c)+\inf_{c\in C}f_B(c)$. 
		
		Using $\{c_n\}_{n=0}^\infty\in IS_{f_A}^P(\alpha,y_0,v_0)$, we get $\displaystyle\sup_{m\in \mathbb{N}}f_A(c_m)< \infty$,\\
		because $\displaystyle\lim_{n\to \infty}f_A(c_n)=\inf_{c\in C}f_A(c)$.
		
		From lemma \ref{l1} it follows, $\displaystyle\sup_{n\in \mathbb{N}}f_B(v_n)<\infty$ and $\displaystyle\sup_{n\in \mathbb{N}}\rho(y_{n},\alpha)< \infty$. By the last two inequalities, $\displaystyle\sup_{m\in \mathbb{N}}f_A(c_m)< \infty$ and (\ref{l3e0}) we obtain for every $k\in \mathbb{N}$ there holds
		\begin{equation}\label{l3e1}
			\begin{array}{rl}
				S\leq&\displaystyle\sup_{m\geq k,\ n\geq k}\rho(y_{n+1},T_A(\alpha,c_m))+I\\[12pt]
				=&\displaystyle\sup_{m\geq k,\ n\geq k}(\rho(y_{n+1},T_A(\alpha,c_m))+I)\\[12pt]
				\leq&\displaystyle \sup_{m\geq k,\ n\geq k}(\lambda(f_A(c_m)+\rho(y_{n},\alpha)+f_B(v_{n}))+(1-\lambda)S)\\[12pt]
				\leq&\displaystyle\lambda\left(\sup_{m\geq k}f_A(c_m)+\sup_{n\geq k}\rho(y_{n},\alpha)+\sup_{n\geq k}f_B(v_{n})\right)+(1-\lambda)S,
			\end{array}
		\end{equation} 
		where $S=S_{f_A,f_B}^C(A,B)$.
		
		Using $\{c_n\}_{n=0}^\infty\in IS_{f_A}^P(\alpha,y_0,v_0)$, we have $\displaystyle \lim_{n\to \infty}f_A(c_n)=\inf_{c\in C} f_A(c)$. From lemma \ref{l2}, we get
		$\displaystyle \lim_{n\to \infty}f_B(v_n)=\inf_{c\in C} f_B(c)$ and\\
		$\displaystyle \lim_{n\to \infty}\rho(y_{n},\alpha)= \dist(A,B)$. By the last three equalities, we can see
		\begin{equation}\label{l3e1.1}
			\left|
			\begin{array}{lclcl}
				\displaystyle \lim_{k\to \infty}\sup_{m\geq k}f_A(c_m)&=&\displaystyle\lim_{k\to \infty}f_A(c_k)&=&\displaystyle\inf_{c\in C} f_A(c)\\[12pt]
				\displaystyle \lim_{k\to \infty}\sup_{n\geq k}f_B(v_n)&=&\displaystyle\lim_{k\to \infty}f_B(v_k)&=&\displaystyle\inf_{c\in C} f_B(c)\\[12pt]
				\displaystyle \lim_{k\to \infty}\sup_{n\geq k}\rho(y_{n},\alpha)&=&\displaystyle\lim_{k\to \infty}\rho(y_{k},\alpha)&=&\displaystyle  \dist(A,B).
			\end{array}
			\right.
		\end{equation} 
		Therefore
		\begin{equation}\label{l3e1.2}
			\lim_{k\to \infty}\lambda\left(\sup_{m\geq k}f_A(c_m)+\sup_{n\geq k}\rho(y_{n},\alpha)+\sup_{n\geq k}f_B(v_{n})\right)+(1-\lambda)S=S.
		\end{equation} 
		Using (\ref{l3e1}) and (\ref{l3e1.2}), we obtain
		$$\displaystyle \lim_{k\to \infty}(\sup_{m\geq k,\ n\geq k}\rho(y_{n+1},T_A(\alpha,c_m))+I)=S.$$
		From the last limit, and $S= \dist(A,B)+I$ it follows
		\begin{equation}\label{l3e2} 
			\lim_{k\to \infty}\sup_{m\geq k,\ n\geq k}\rho(y_{n+1},T_A(\alpha,c_m))= \dist(A,B).
		\end{equation} 
		By (\ref{l3e1.3}), (\ref{l3e-1}), and $(T_A,T_B,H_A,H_B,f_A,f_B)\in CEF_C^P(A,B)$, we get 
		$$
		\begin{array}{rl}
			S=&\displaystyle  \dist(A,B)+\inf_{c\in C}f_A(c)+\inf_{c\in C}f_B(c)\\
			\leq&\displaystyle f_A(H_A(\alpha,c_n))+  \dist(A,B)+\inf_{c\in C}f_B(c)\\
			\leq&\displaystyle f_A(H_A(\alpha,c_n))+\rho(y_{n+1},T_A(\alpha,c_n))+f_B(v_{n+1})\\
			\leq&\displaystyle\lambda(f_A(c_n)+\rho(y_{n},\alpha)+f_B(v_{n}))+(1-\lambda)S.
		\end{array}
		$$
		Using the last chain of inequalities, and (\ref{l3e1.1}), we can observe
		\begin{equation}\label{l3e2.1}  
			\begin{array}{lll}
				\displaystyle\lim_{n\to \infty}f_A(c_n)&=&\displaystyle\inf_{c\in C} f_A(c)\\[12pt]
				\displaystyle\lim_{n\to \infty}f_B(v_n)&=&\displaystyle\inf_{c\in C} f_B(c)\\[12pt]
				\displaystyle\lim_{n\to \infty}f_A(H_A(\alpha,c_n))&=&\displaystyle\inf_{c\in C} f_A(c).\\
			\end{array}
		\end{equation} 
		From lemma \ref{l2}, we get $\displaystyle\lim_{n\to \infty} \rho(y_{n},\alpha)=\lim_{n\to \infty} \rho(y_{n+1},\alpha)= \dist(A,B)$, i.e.,
		\begin{equation}\label{l3e3}  
			\lim_{k\to \infty}\sup_{m\geq k,\ n\geq k}\rho(y_{n+1},\alpha)=\lim_{k\to \infty}\sup_{n\geq k}\rho(y_{n+1},\alpha)= \dist(A,B).
		\end{equation} 
		Using definition \ref{contraction}:\ref{C-1} and $\{c_n\}_{n=0}^\infty\in IS_{f_A}^P(\alpha,y_0,v_0)$,\\
		it follows $(\alpha,y_{m+1},c_n,v_{m+1})\in P$ and \\ $(T_A(\alpha,c_n),y_{m+1},H_A(\alpha,c_n),v_{m+1})\in P$ for every $n,m\in \mathbb{N}$. By the last two inclusions, (\ref{l3e2}), (\ref{l3e2.1}), (\ref{l3e3}), lemma \ref{l0.1}, and $(A,B)\in CD_X^P(f_A,f_B)$, we can observe that
		$$
		\alpha=\lim_{n\to \infty}\alpha=\lim_{n\to \infty}T_A(\alpha,c_n).
		$$
	\end{proof}
	\begin{lem}\label{l4}
		Let $(X,\rho)$ be a metric space, $A$ and $B$ be subsets of $X$, $C$ be an arbitrary set, $P\subseteq A\times B\times C^2$, and $(T_A,T_B,H_A,H_B,f_A,f_B)\in CEF_C^P(A,B)$. Let $\alpha\in A$ satisfy $\displaystyle\lim_{n\to \infty}T_A(\alpha,c_n)=\alpha$ where\\
		$\{c_n\}_{n=0}^\infty\in IS_{f_A}^P(\alpha,y_0,v_0)$, and $(y_0,v_0)\in B\times C$. Let $\{y_n\}_{n=0}^\infty$ and $\{v_n\}_{n=0}^\infty$ be the iterated sequences of $(T_B,H_B)$ with an initial guess $(y_0,v_0)$. Then 
		$$\lim_{n\to \infty}\rho(y_n,\alpha)= \dist(A,B)$$
	\end{lem}
	\begin{proof}
		From definition \ref{contraction}:\ref{C-1}, definition \ref{def4} and $\{c_n\}_{n=0}^\infty\in IS_{f_A}^P(\alpha,y_0,v_0)$, we get for every $n,m\in \mathbb{N}$ there holds
		$$
		(\alpha,y_m,c_n,v_m)\in P.
		$$
		By the last inclusion, and $(T_A,T_B,H_A,H_B,f_A,f_B)\in CEF_C^P(A,B)$, we can see that for every $n\in \mathbb{N}$, the chain of inequalities
		$$
		\begin{array}{rl}
			&\rho(y_{n+1},\alpha)-\rho(\alpha,T_A(\alpha,c_n))+I\\
			\leq&\rho(y_{n+1},T_A(\alpha,c_n))+I\\
			\leq&f_A(H_A(\alpha,c_n))+\rho(y_{n+1},T_A(\alpha,c_n))+f_B(v_{n+1})\\
			\leq&\lambda(f_A(c_n)+\rho(y_{n},\alpha)+f_B(v_{n}))+(1-\lambda)S_{f_A,f_B}^C(A,B)\\
			=&\lambda(\rho(y_{n},\alpha)+f_B(v_{n})+f_A(c_n)-I)+(1-\lambda) \dist(A,B)+I,
		\end{array}
		$$
		is true, where $\displaystyle I=\inf_{c\in C}f_A(c)+\inf_{c\in C}f_B(c)$. Therefore
		\begin{equation}\label{l4e1}
			\rho(y_{n+1},\alpha)\leq\lambda\rho(y_{n},\alpha)+Q_n+(1-\lambda) \dist(A,B),
		\end{equation}
		where $\displaystyle Q_n=\lambda(f_B(v_{n})+f_A(c_n)-I)+\rho(\alpha,T_A(\alpha,c_n))$.\\
		Using $(\alpha,y_0,c_0,v_0)\in P$, lemma \ref{l2}, $\{c_n\}_{n=0}^\infty\in IS_{f_A}^P(\alpha,y_0,v_0)$, we can observe that $\displaystyle\lim_{n\to \infty}f_A(c_n)+\lim_{n\to \infty}f_B(v_{n})=I$. From the last equality, the continuity of the metric function, and $\displaystyle\lim_{n\to \infty}T_A(\alpha,c_n)=\alpha$, we obtain $\displaystyle\lim_{n\to\infty}Q_n=0$. Thus
		\begin{equation}\label{l4e2}
			\lim_{n\to\infty}J_n=0,
		\end{equation}
		where $\displaystyle J_n=\sup_{k\geq n}Q_k$.\\
		Applying (\ref{l4e1}) $m$ times consecutively, we get
		\begin{equation}\label{l4e2.1}
			\begin{array}{rl}
				W_8=&\rho(y_{n+m},\alpha)\\
				\leq&\lambda\rho(y_{n+m-1},\alpha)+Q_{n+m-1}+(1-\lambda) \dist(A,B)\\
				\leq&\lambda^2\rho(y_{n+m-2},\alpha)+  Q_{n+m-1}+\lambda Q_{n+m-2}+(1-\lambda^2) \dist(A,B)\\
				&\dots\\
				\leq&\displaystyle\lambda^{m-2}\rho(y_{n+2},\alpha)+\sum_{k=1}^{m-2}\lambda^{k-1}Q_{n+m-k}+(1-\lambda^{m-2}) \dist(A,B)\\[16pt]
				\leq&\displaystyle\lambda^{m-1}\rho(y_{n+1},\alpha)+\sum_{k=1}^{m-1}\lambda^{k-1}Q_{n+m-k}+(1-\lambda^{m-1}) \dist(A,B)\\[16pt]
				\leq&\displaystyle\lambda^{m}\rho(y_{n},\alpha)+\sum_{k=1}^{m}\lambda^{k-1}Q_{n+m-k}+(1-\lambda^{m}) \dist(A,B).\\
			\end{array}
		\end{equation}
		
		By the definition of $J_n$, we can see for every $n\in \mathbb{N}$ and $k\in \mathbb{N}$ so that $n\leq k$, there holds $Q_k\leq J_n$. From the last inequality and (\ref{l4e2.1}), it follows for every $n\in \mathbb{N}$ and $m\in \mathbb{N}$ the inequality
		\begin{equation}\label{l4e3}
			\rho(y_{n+m},\alpha)\leq \lambda^{m}\rho(y_{n},\alpha)+ \frac{J_n}{1-\lambda}+(1-\lambda^{m}) \dist(A,B),
		\end{equation}
		is valid. For any $n\in \mathbb{N}$, the limit
		$$\lim_{m\to \infty}\left(\lambda^{m}\rho(y_{n},\alpha)+(1-\lambda^{m}) \dist(A,B)\right)= \dist(A,B),$$
		is true. By the last equality and (\ref{l4e2}) we obtain, for every $\varepsilon>0$ and $\delta>0$, we can choose $N_\delta\in \mathbb{N}$ and $M(N_\delta,\varepsilon)\in \mathbb{N}$, such that if $m\geq M(N_\delta,\varepsilon)$ then
		$$\left|
		\begin{array}{l}
			\displaystyle\frac{J_{N_\delta}}{1-\lambda}<\delta\\[12pt]
			\displaystyle\lambda^{m}\rho(y_{N_\delta},\alpha)+(1-\lambda^{m}) \dist(A,B)< \dist(A,B)+\varepsilon.
		\end{array}
		\right.$$
		Using the last system of inequalities, (\ref{l4e3}), $\alpha\in A$ and $\{y_n\}_{n=0}^\infty\subseteq B$, we get for any $m\geq M(N_\delta,\varepsilon)$, there holds
		$$ \dist(A,B)\leq\rho(y_{(N_\delta+m)},\alpha)\leq  \dist(A,B)+\varepsilon+\delta.$$ 
		Consequently
		$$\lim_{n\to \infty}\rho(y_n,\alpha)= \dist(A,B).$$
	\end{proof}
	\begin{lem}\label{l5}
		Let $(X,\rho)$ be a metric space, $A$ and $B$ be subsets of $X$, $C$ be an arbitrary set, $P\subseteq A\times B\times C^2$, $(A,B)\in CD_X^P(f_A,f_B)$, and\\
		$(T_A,T_B,H_A,H_B,f_A,f_B)\in CEF_C^P(A,B)$. Let $(x_0,y_0,u_0,v_0)\in P$, $\{x_n\}_{n=0}^\infty$ and $\{u_n\}_{n=0}^\infty$ be the iterated sequences of $(T_A,H_A)$ with an initial guess $(x_0,u_0)$, and $\displaystyle\lim_{n\to\infty}x_n=\alpha$. Then there is no $\beta\in A$, so that\\
		$\displaystyle\lim_{n\to \infty}T_A(\beta,c_n)=\beta$ and $\beta\neq \alpha$ , where $\{c_n\}_{n=0}^\infty\in IS_{f_A}^P(\beta,y_0,v_0)$.
	\end{lem}
	\begin{proof}
		We will prove this lemma through the assumption that there exists such $\beta$, and the showing of a contradiction.
		
		Let $\beta\in A$ be such that $\beta\neq \alpha$ and $\displaystyle\lim_{n\to \infty}T_A(\beta,c_n)=\beta$, where\\
		$\{c_n\}_{n=0}^\infty\in IS_{f_A}^P(\beta,y_0,v_0)$.
		Let $\{y_n\}_{n=0}^\infty$ and $\{v_n\}_{n=0}^\infty$ be the iterated sequences of\\
		$(T_B,H_B)$ with an initial guess $(y_0,v_0)$.
		
		From $\{c_n\}_{n=0}^\infty\in IS_{f_A}^P(\beta,y_0,v_0)$, lemma \ref{l2} and lemma \ref{l4} we get
		\begin{equation}\label{l5e1}
			\left|
			\begin{array}{l}
				\displaystyle \lim_{k\to \infty}\sup_{n\geq k,\ m\geq k}\rho(x_n,y_m)= \dist(A,B)\\
				\displaystyle  \dist(A,B)=\lim_{n\to \infty}\rho(\beta,y_{n})=\lim_{k\to \infty}\sup_{n\geq k,\ m\geq k}\rho(\beta,y_m).
			\end{array}
			\right.
		\end{equation}
		Using lemma \ref{l2} and $\{c_n\}_{n=0}^\infty\in IS_{f_A}^P(\beta,y_0,v_0)$, it follows
		\begin{equation}\label{l5e2}
			\begin{array}{rcl}
				\displaystyle\lim_{n\to \infty}f_A(c_n)&=&\displaystyle\inf_{c\in C} f_A(c)\\
				\displaystyle\lim_{n\to \infty}f_A(u_n)&=&\displaystyle\inf_{c\in C} f_A(c)\\
				\displaystyle\lim_{n\to \infty}f_B(v_n)&=&\displaystyle\inf_{c\in C} f_B(c).
			\end{array}
		\end{equation}
		By definition \ref{contraction}:\ref{C-1} and $\{c_n\}_{n=0}^\infty\in IS_{f_A}^P(\beta,y_0,v_0)$, we obtain\\
		$(x_n,y_m,u_n,v_m)\in P$ and $(\beta,y_m,c_n,v_m)\in P$ for every $n,m\in \mathbb{N}$.
		From the last two inclusions, lemma \ref{l0.1}, (\ref{l5e1}), (\ref{l5e2}), and $(A,B)\in CD_X^P(f_A,f_B)$, we can observe that
		$$
		\alpha=\lim_{n\to \infty}x_n=\lim_{n\to \infty}\beta=\beta.
		$$
		This is a contradiction to $\beta\neq \alpha$
		
	\end{proof}﻿
	
	\section{\normalsize Proof of theorem \ref{t1}}
	
	\begin{proof}
		From lemma \ref{l2}, it follows for every $(x_0,y_0,u_0,v_0)\in P$ there exists $\displaystyle\lim_{n\to\infty}x_n=\alpha\in A$ and the equalities
		$$
		\begin{array}{l}
			\displaystyle\lim_{n\to \infty}\rho(y_n,\alpha)= \dist(A,B)\\
			\displaystyle\lim_{n\to \infty} f_A(u_n)= \inf_{c\in C} f_A(c)\\
			\displaystyle\lim_{n\to \infty}f_B(v_n)= \inf_{c\in C} f_B(c),
		\end{array}
		$$
		are true. Where $\{x_n\}_{n=0}^\infty$ and $\{u_n\}_{n=0}^\infty$ are the iterated sequences of $(T_A,H_A)$ with an initial guess $(x_0,u_0)$,  $\{y_n\}_{n=0}^\infty$ and $\{v_n\}_{n=0}^\infty$ are the iterated sequences of $(T_B,H_B)$ with an initial guess $(y_0,v_0)$.
		
		By lemma \ref{l2.1}, we get if $(x_0,y_0,u_0.v_0)\in P$ and\\
		$(z_0,y_0,t_0,v_0)\in P$, then $\displaystyle\lim_{n\to \infty}x_n=\lim_{n\to \infty}z_n\in A$ where $\{x_n\}_{n=0}^\infty$ and $\{u_n\}_{n=0}^\infty$ are the iterated sequences of $(T_A,H_A)$ with an initial guess $(x_0,u_0)$, and $\{z_n\}_{n=0}^\infty$ and $\{t_n\}_{n=0}^\infty$ are the iterated sequences of $(T_A,H_A)$ with an initial guess $(z_0,t_0)$.
		
		Using lemma \ref{l3} we obtain if $(x_0,y_0,u_0,v_0)\in P$, $\{x_n\}_{n=0}^\infty$ and $\{u_n\}_{n=0}^\infty$ are the iterated sequences of $(T_A,H_A)$ with an initial guess $(x_0,u_0)$, $\alpha=\displaystyle\lim_{n\to \infty}x_{n}$ and $\{c_n\}_{n=0}^\infty\in IS_{f_A}^P(\alpha,y_0,v_0)$, then $\displaystyle\lim_{n\to \infty}T_A(\alpha,c_n)=\alpha$.
		
		From lemma \ref{l5}, we can observe if $(x_0,y_0,u_0,v_0)\in P$,$\{x_n\}_{n=0}^\infty$ and $\{u_n\}_{n=0}^\infty$ are the iterated sequences of $(T_A,H_A)$ with an initial guess $(x_0,u_0)$, and $\displaystyle\lim_{n\to\infty}x_n=\alpha$, then there is no $\beta\in A$ so that $\displaystyle\lim_{n\to \infty}T_A(\beta,c_n)=\beta$ and $\beta\neq \alpha$ , where $\{c_n\}_{n=0}^\infty\in IS_{f_A}^P(\beta,y_0,v_0)$.
	\end{proof}
	
	\section{\normalsize Examples and Applications}
	
	\subsection{Examples}
	
	We will now introduce a more restrictive version of the $CD$ property, which is, however easier to prove.
	
	\begin{defn}
		Let $(X,\rho)$ be a metric space, $A$ and $B$ be subsets of $X$, and for every pair of sequences $\{x_n\}_{n=1}^{\infty}\in A$ and $\{y_n\}_{n=1}^{\infty}\in B$, such that $\displaystyle\lim_{k\to \infty}\sup_{n\geq k,\ m\geq k}\rho(x_n,y_m)= \dist(A,B)$, there holds $\{x_n\}_{n=1}^{\infty}$ converges in $A$. Then we say that the ordered pair $(A,B)$ satisfies the property Convergent by Distance $(CD)$.
	\end{defn}
	
	We will continue with a theorem that shows one of the possible relationships between the $UC$ and $CD$ properties.
	
	\begin{thm}\label{t0}
		Let $(X,\rho)$ be a metric space, $A$ and $B$ be subsets of $X$, $(A,\rho)$ be complete, and the ordered pair $(A,B)$ satisfy the property $UC$. Then the ordered pair $(A,B)$ satisfies the property $CD$.
	\end{thm}
	\begin{proof}
		We will prove this theorem through the assumption that the ordered pair $(A,B)$ does not satisfy the property $CD$, and the showing of a contradiction.
		
		Let us assume that the ordered pair $(A,B)$ does not satisfy the property $CD$. Then there are sequences $\{x_n\}_{n=1}^\infty\subseteq A$ and $\{y_n\}_{n=1}^\infty\subseteq B$ so that
		\begin{equation}\label{t0e1}
			\lim_{k\to \infty}\sup_{n\geq k,\ m\geq k}\rho(x_n,y_m)= \dist(A,B),
		\end{equation} 
		and $\{x_n\}_{n=1}^\infty$ does not converge in $A$. By assumption, $(A,\rho)$ is complete. Therefore, there exists $\varepsilon>0$ such that for every $N\in \mathbb{N}$, there are $n_N>m_N\geq N$ so that $\rho(x_{n_N},x_{m_N})>\varepsilon$, i.e.,
		\begin{equation}\label{t0e2}
			\lim_{N\to \infty}\rho(x_{n_N},x_{m_N})\neq 0.
		\end{equation} 
		From $\{x_n\}_{n=1}^\infty\subseteq A$, $\{y_n\}_{n=1}^\infty\subseteq B$, and $n_N>m_N\geq N$, it follows that for every $N\in \mathbb{N}$ there holds 
		$$
		\begin{array}{l}
			\displaystyle  \dist(A,B)\leq \rho(x_{n_N},y_{n_N})\leq \sup_{p\geq N,\ q\geq N}\rho(x_p,y_q)\\
			\displaystyle  \dist(A,B)\leq \rho(x_{m_N},y_{n_N})\leq \sup_{p\geq N,\ q\geq N}\rho(x_p,y_q).
		\end{array}
		$$
		Using the last two inequalities and (\ref{t0e1}), we obtain 
		$$
		\begin{array}{l}
			\displaystyle \lim_{N\to \infty}\rho(x_{n_N},y_{n_N})=  \dist(A,B)\\
			\displaystyle \lim_{N\to \infty}\rho(x_{m_N},y_{n_N})=  \dist(A,B).
		\end{array}
		$$
		the last two limits, and (\ref{t0e2}),
		contradicts the assumption that the ordered pair $(A,B)$ satisfies the property $UC$.
	\end{proof}
	
	From theorem \ref{t0} it follows that if we add the condition that $A$ is a complete subspace to the points of example \ref{e0} they become examples for ordered pairs satisfying the $CD$ property.
	
	In the next example we will use the convention $2^{\lfloor\log_{2}x\rfloor} = 0$ for $x=0$.
	
	\begin{ex}\label{e1}
		Let us consider the metric space $(\mathbb{R},\rho(\cdot,\cdot))$, where $\rho(\cdot,\cdot)=|\cdot-\cdot|$, let $A=[0;\infty)$, $B=(-\infty,-1]$, $C=A\cup B$, and $P=\{(a,b,a,b)\in A\times B \times C^2:a\in A,\ b\in B\}$. Let us define
		$$\alpha(x)=\left\{
		\begin{array}{ll}
			\lfloor\log_{2}x\rfloor\ \mod\ 2&:x> 0\\
			0&:x= 0
		\end{array}
		\right.$$
		$T:A\to A$ as $\displaystyle Tx=2x\alpha(x)+\frac{1}{4}(x-2^{\lfloor\log_{2}x\rfloor})(1-\alpha(x))$\\
		$T_A:A\times C\to A$ as $T_A(a,c)=Ta$ \\
		$T_B:B\times C\to B$ as $\displaystyle T_B(b,c)=\frac{b+1}{8}+\frac{15}{8}(b+1)\alpha(-b-1)-1$\\
		$H_A:A\times C\to C$ as $H_A(a,c)=Ta$\\
		$H_B:B\times C\to C$ as $H_B(b,c)=T_B(b,c)$\\
		$f_A: C\to \mathbb{R}$ as 
		$$
		f_A(c)=\left\{
		\begin{array}{ll}
			4c\alpha(c)&:c\geq 0\\
			0&:c\leq -1
		\end{array}
		\right.
		$$
		$f_B: C\to \mathbb{R}$ as 
		$$
		f_B(c)=\left\{
		\begin{array}{ll}
			0&:c\geq 0\\
			-4(c+1)\alpha(-c-1)&:c\leq -1.
		\end{array}
		\right.
		$$
		Then $(T_A,T_B,H_A,H_B,f_A,f_B)\in CEF_C^P(A,B)$, the sequence $\{T^na\}_{n=0}^{\infty}$ converges to $0$, for every $a\in A$, and $x=0$ is a unique fixed point of $T$ in $\mathbb{L}=\{x\in \mathbb{R}:\alpha(x)=0,\ x\geq 0\}$. 
	\end{ex}
	\begin{proof}{\rm
			First, let us consider some properties of $\alpha$
			\begin{equation}\label{e1e0.1}
				\alpha(x)\in \{0,1\}\leq1\ for\ every\ x\in \mathbb{R}.
			\end{equation}
			For every $x>0$ the equality
			\begin{equation}\label{e1e0.2}
				\alpha(2x)=\lfloor\log_{2}2x\rfloor\ \mod\ 2=(1+\lfloor\log_{2}x\rfloor)\ \mod\ 2=1-\alpha(x),
			\end{equation}
			is true. From (\ref{e1e0.1}), it follows that $\displaystyle \inf_{c\in C}f_A(c)=0. $ By similar arguments, we can observe that $\displaystyle \inf_{c\in C}f_B(c)=0$.
			
			Now let us consider $x-2^{\lfloor\log_{2}x\rfloor}$. The inequality
			\begin{equation}\label{e1e0.3}
				x-2^{\lfloor\log_{2}x\rfloor}\leq x-2^{\log_{2}x-1}= x-\frac{1}{2}x=\frac{1}{2}x,
			\end{equation}
			is valid.  We can see that $ \dist(A,B)=1$.
			
			It is easy to check that, $A$ is convex. Then, by example corollary \ref{c-0}, it follows that the ordered pair $(A,B)$ satisfies the $UC$ property. We can observe that $A$ is closed. Then, using theorem \ref{t0}, we get that the ordered pair $(A,B)$ satisfies the $CD$ property. Therefore, $(A,B)\in CD_\mathbb{R}^P(f_A,f_B)$.
			
			We can see that the condition from definition \ref{contraction}:\ref{C-1} holds true.
			
			Let us Consider
			\begin{equation}\label{e1e1}
				\rho(T_A(x,u),T_B(y,v))+f_A(H_A(x,u))+f_B(H_B(y,v)).
			\end{equation}
			For an arbitrary $(x,y,u,v)\in P$, there holds $u=x$ and $v=y$. Thus there are four cases:
			\begin{enumerate}[label=Case-\arabic*:]
				\item $\alpha(x)=0$ and $\alpha(-1-y)=0$. Then from $x\geq0$, $y<0$, $x-2^{\lfloor\log_{2}x\rfloor}\geq0$, $\displaystyle-\frac{y+1}{8}\geq 0$, (\ref{e1e0.1}), and (\ref{e1e0.3}), it follows (\ref{e1e1}) would look like
				$$
				\begin{array}{rl}
					&\displaystyle\rho\left(T_A\left(x,u\right),T_B\left(y,v\right)\right)+f_A\left(H_A\left(x,u\right)\right)+f_B\left(H_B\left(y,v\right)\right)\\[16pt]
					=&\displaystyle\left|\frac{1}{4}\left(x-2^{\lfloor\log_{2}x\rfloor}\right)-\frac{y+1}{8}+1\right|+\left(x-2^{\lfloor\log_{2}x\rfloor}\right)\alpha\left(\frac{1}{4}\left(x-2^{\lfloor\log_{2}x\rfloor}\right)\right)\\
					&\displaystyle-\frac{y+1}{2}\alpha\left(-\frac{y+1}{8}\right)\\[16pt]
					\leq&\displaystyle\frac{1}{4}\left(x-2^{\lfloor\log_{2}x\rfloor}\right)-\frac{y+1}{8}+1+\left(x-2^{\lfloor\log_{2}x\rfloor}\right)-\frac{y+1}{2}\\[16pt]
					=&\displaystyle\frac{5}{4}\left(x-2^{\lfloor\log_{2}x\rfloor}\right)-\frac{5}{8}y+\left(1-\frac{5}{8}\right)\\[16pt]
					\leq&\displaystyle\frac{5}{8}\left(x-y\right)+\left(1-\frac{5}{8}\right).1\\[16pt]
					\leq&\displaystyle\frac{5}{8}\left(|x-y|+0+0\right)+\left(1-\frac{5}{8}\right) \dist\left(A,B\right)\\[16pt]
					=&\displaystyle\frac{5}{8}\left(\rho\left(x,y\right)+f_A\left(u\right)+f_B\left(v\right)\right)+\left(1-\frac{5}{8}\right) \dist(A,B)
				\end{array}
				$$
				\item $\alpha(x)=1$ and $\alpha(-1-y)=0$. Then using $x\geq0$, $y<0$, $\displaystyle-\frac{y+1}{8}\geq 0$, (\ref{e1e0.1}), and (\ref{e1e0.2}), we get (\ref{e1e1}) would look like
				$$
				\begin{array}{rl}
					&\displaystyle\rho\left(T_A\left(x,u\right),T_B\left(y,v\right)\right)+f_A\left(H_A\left(x,u\right)\right)+f_B\left(H_B\left(y,v\right)\right)\\[16pt]
					=&\displaystyle\left|2x-\frac{y+1}{8}+1\right|+2x\alpha(2x)-\frac{y+1}{2}\alpha\left(-\frac{y+1}{8}\right)\\[16pt]
					\leq&\displaystyle2x-\frac{y+1}{8}+1-\frac{y+1}{2}\\[16pt]
					=&\displaystyle\frac{5}{8}x-\frac{5}{8}y+\frac{11}{8}x+\left(1-\frac{5}{8}\right)\\[16pt]
					\leq&\displaystyle\frac{5}{8}\left(x-y+4x\right)+\left(1-\frac{5}{8}\right).1\\[16pt]
					=&\displaystyle\frac{5}{8}\left(|x-y|+4x+0\right)+\left(1-\frac{5}{8}\right) \dist\left(A,B\right)\\[16pt]
					=&\displaystyle\frac{5}{8}\left(\rho\left(x,y\right)+f_A\left(u\right)+f_B\left(v\right)\right)+\left(1-\frac{5}{8}\right) \dist(A,B)
				\end{array}
				$$
				\item  $\alpha(x)=0$ and $\alpha(-1-y)=1$. Then from $x\geq0$, $y<0$, $x-2^{\lfloor\log_{2}x\rfloor}\geq0$, $\displaystyle-2(y+1)\geq 0$, (\ref{e1e0.1}), (\ref{e1e0.2}), and (\ref{e1e0.3}), we obtain (\ref{e1e1}) would look like
				$$
				\begin{array}{rl}
					&\displaystyle\rho\left(T_A\left(x,u\right),T_B\left(y,v\right)\right)+f_A\left(H_A\left(x,u\right)\right)+f_B\left(H_B\left(y,v\right)\right)\\[16pt]
					=&\displaystyle\left|\frac{1}{4}\left(x-2^{\lfloor\log_{2}x\rfloor}\right)-2(y+1)+1\right|\\
					&\displaystyle+\left(x-2^{\lfloor\log_{2}x\rfloor}\right)\alpha\left(\frac{1}{4}\left(x-2^{\lfloor\log_{2}x\rfloor}\right)\right)-2(y+1)\alpha\left(2(-y-1)\right)\\[16pt]
					\leq&\displaystyle\frac{1}{4}\left(x-2^{\lfloor\log_{2}x\rfloor}\right)-2(y+1)+1+\left(x-2^{\lfloor\log_{2}x\rfloor}\right)\\[16pt]
					=&\displaystyle\frac{5}{4}\left(x-2^{\lfloor\log_{2}x\rfloor}\right)-\frac{5}{8}y-\frac{11}{8}(y+1)+\left(1-\frac{5}{8}\right)\\[16pt]
					\leq&\displaystyle\frac{5}{8}\left(x-y-4(y+1)\right)+\left(1-\frac{5}{8}\right).1\\[16pt]
					\leq&\displaystyle\frac{5}{8}\left(|x-y|+0-4(y+1)\right)+\left(1-\frac{5}{8}\right) \dist\left(A,B\right)\\[16pt]
					=&\displaystyle\frac{5}{8}\left(\rho\left(x,y\right)+f_A\left(u\right)+f_B\left(v\right)\right)+\left(1-\frac{5}{8}\right) \dist(A,B)
				\end{array}
				$$
				\item  $\alpha(x)=1$ and $\alpha(-1-y)=1$. Then by $x\geq0$, $y<0$, $\displaystyle-2(y+1)\geq 0$, and (\ref{e1e0.2}), we can observe (\ref{e1e1}) would look like
				$$
				\begin{array}{rl}
					&\displaystyle\rho\left(T_A\left(x,u\right),T_B\left(y,v\right)\right)+f_A\left(H_A\left(x,u\right)\right)+f_B\left(H_B\left(y,v\right)\right)\\[16pt]
					=&\displaystyle\left|2x-2(y+1)+1\right|+2x\alpha(2x)-2(y+1)\alpha\left(2(-y-1)\right)\\[16pt]
					=&\displaystyle2x-2(y+1)+1\\[16pt]
					=&\displaystyle\frac{5}{8}x-\frac{5}{8}y+\frac{11}{8}x-\frac{11}{8}(y+1)+\left(1-\frac{5}{8}\right)\\[16pt]
					\leq&\displaystyle\frac{5}{8}\left(|x-y|+4x-4(y+1)\right)+\left(1-\frac{5}{8}\right) \dist\left(A,B\right)\\[16pt]
					=&\displaystyle\frac{5}{8}\left(\rho\left(x,y\right)+f_A\left(u\right)+f_B\left(v\right)\right)+\left(1-\frac{5}{8}\right) \dist(A,B).
				\end{array}
				$$
			\end{enumerate}
			In all cases, we get 
			\begin{equation}\label{e1e1.1}
				\begin{array}{rcl}
					W_9&=&\rho\left(T_A\left(x,u\right),T_B\left(y,v\right)\right)+f_A\left(H_A\left(x,u\right)\right)+f_B\left(H_B\left(y,v\right)\right)\\
					&\leq&\displaystyle\frac{5}{8}\left(\rho\left(x,y\right)+f_A\left(u\right)+f_B\left(v\right)\right)+\left(1-\frac{5}{8}\right) \dist(A,B).
				\end{array}
			\end{equation}
			
			It is easy to check that \ref{contraction}:\ref{C-1} holds, the inequality (\ref{e1e1.1}), is true, $\displaystyle \inf_{c\in C}f_A(c)=0$, and $\displaystyle \inf_{c\in C}f_B(c)=0$. Consequently, $(T_A,T_B,H_A,H_B,f_A,f_B)\in CEF_C^P(A,B)$, $(A,B)\in CD_X^P(f_A,f_B)$. Thus, we can apply theorem \ref{t1}.
			
			For every $a\in A$, there holds $(a,-1,a,-1)\in P$, $(0,-1,0,-1)\in P$, 
			$\{T^na\}_{n=0}^{\infty}$ and $\{T^na\}_{n=0}^{\infty}$ are the iterated sequences of $(T_A,H_A)$ for an arbitrary initial guess $(a,a)\in A\times C$, and $\{0\}_{n=0}^{\infty}$, and $\{0\}_{n=0}^{\infty}$ are the iterated sequences of $(T_A,H_A)$ for an initial guess $(0,0)$. Then, from theorem \ref{t1} consequence \ref{t1-cons2}, it follows that for every $a\in A$ the limit\\
			$\ds \lim_{n\to \infty}T^na=\lim_{n\to \infty}0=0$, is true.
			
			Let us suppose that there is $\beta\in\mathbb{L}$ such that $\beta=T\beta$ and $\beta\neq0$. By $\beta\in\mathbb{L}$, we get 
			\begin{equation}\label{e1e2}
				\displaystyle\lim_{n\to \infty}f_A(\beta)=4\beta\alpha(\beta)=0=\inf_{c\in C}f_A(c).
			\end{equation}
			It is easy to see that $(\beta,-1,\beta,-1)\in P$. Then, using \ref{e1e2}, it follows $\{\beta\}_{n=0}^{\infty}\in IS_{f_A}^P(\beta,-1,-1)$. From $T\beta=\beta$, we obtain $\displaystyle\lim_{n\to \infty}T_A(\beta,\beta)=\beta$. We have $\{0\}_{n=0}^{\infty}$, $\{0\}_{n=0}^{\infty}$ are the iterated sequences of $(T_A,H_A)$ for an initial guess $(0,0)$, and $\displaystyle\lim_{n\to \infty}0=0$. Then, $\{\beta\}_{n=0}^{\infty}\in IS_{f_A}^P(\beta,-1,-1)$, and $\displaystyle\lim_{n\to \infty}T_A(\beta,\beta)=\beta$, contradict theorem \ref{t1} consequence \ref{t1-cons4}. Consequently,
			$x=0$ is a unique fixed point of $T$ in $\mathbb{L}$. }
	\end{proof}
	
	\subsection{Application} 
	
	The following proposition is in order to facilitate the establishment of the $CD$ property in certain cases of a Cartesian product.
	
	\begin{prop}\label{t0.1}
		Let $(X_1,\rho_1)$ and $(X_2,\rho_2)$ be metric spaces, $A_1,B_1\subseteq X_1$,  $A_2,B_2\subseteq X_2$, and the ordered pairs $(A_1,B_1)$ and $(A_2,B_2)$ satisfy the property $CD$. Then the ordered pair $(A_1\times A_2, B_1\times B_2)$ satisfies the property $CD$ in the space $(X_1\times X_2, d(\cdot,\cdot))$, where $d=\rho_1+\rho_2$.
	\end{prop}
	\begin{proof}
		Let $\{p_n\}_{n=0}^{\infty}=\{(a_n,b_n)\}_{n=0}^{\infty}\subseteq A_1\times A_2$, $\{q_n\}_{n=0}^{\infty}=\{(c_n,d_n)\}_{n=0}^{\infty}\subseteq B_1\times B_2$, and 
		$$
		\lim_{k\to \infty}\sup_{n\geq k,\ m\geq k}d(p_n,q_m)= \dist(A_1\times A_2,B_1\times B_2).
		$$
		Therefore 
		\begin{equation}\label{t0.1e1}
			\lim_{k\to \infty}\sup_{n\geq k,\ m\geq k}(\rho_1(a_n,c_m)+\rho_2(b_n,d_m))= \dist(A_1,B_1)+ \dist(A_2,B_2).
		\end{equation} 
		
		Then, from $\{a_n\}_{n=0}^{\infty}\subseteq A_1$, $\{b_n\}_{n=0}^{\infty}\subseteq A_2$, $\{c_n\}_{n=0}^{\infty}\subseteq B_1$, and \\
		$\{d_n\}_{n=0}^{\infty}\subseteq B_2$, it follows that for every $n,m\in \mathbb{N}_0$, the inequalities\\
		$\rho_1(a_n,c_m)\geq \dist(A_1,B_1)$ and $\rho_2(b_n,d_m)\geq \dist(A_2,B_2)$, are true. Thus, if $\displaystyle\inf_{n,m\in\mathbb{N}}\rho_1(a_n,c_m)>\dist(A_1,B_1)$, or $\displaystyle\inf_{n,m\in\mathbb{N}}\rho_2(b_n,d_m)>\dist(A_2,B_2)$, then there is $\varepsilon>0$, so that for any $k\in \mathbb{N}_0$ there holds
		$$
		\begin{array}{l}
			\ds\sup_{n\geq k,\ m\geq k}(\rho_1(a_n,c_m)+\rho_2(b_n,d_m))\geq\inf_{n,m\in\mathbb{N}}\rho_1(a_n,c_m)+\inf_{n,m\in\mathbb{N}}\rho_2(b_n,d_m)\\[16pt]
			\geq\dist(A_1,B_1)+ \dist(A_2,B_2)+\varepsilon.
		\end{array}
		$$
		The last inequality contradicts (\ref{t0.1e1}). Consequently
		\begin{equation}\label{t0.1e2}
			\begin{array}{l}
				\displaystyle\inf_{n,m\in\mathbb{N}}\rho_1(a_n,c_m)= \dist(A_1,B_1)\\
				\displaystyle\inf_{n,m\in\mathbb{N}}\rho_2(b_n,d_m)= \dist(A_2,B_2).
			\end{array}
		\end{equation} 
		
		By (\ref{t0.1e1}), (\ref{t0.1e2}) and lemma \ref{l0.0}, we get
		$$
		\left|
		\begin{array}{l}
			\displaystyle\lim_{k\to \infty}\sup_{n\geq k,\ m\geq k}\rho_1(a_n,c_m)= \dist(A_1,B_1)\\
			\displaystyle \lim_{k\to \infty}\sup_{n\geq k,\ m\geq k}\rho_2(b_n,d_m)= \dist(A_2,B_2).
		\end{array}
		\right.
		$$
		Using the last system, and the assumption that the ordered pairs $(A_1,B_1)$ and $(A_2,B_2)$ satisfy the property $CD$, we obtain
		$$
		\begin{array}{l}
			\displaystyle\lim_{n\to \infty}a_n=\mu \in A_1\\
			\displaystyle\lim_{n\to \infty}b_n=\nu \in A_2.
		\end{array}
		$$
		Therefore
		$$
		\lim_{n\to \infty}d(p_n,(\mu,\nu))=\lim_{n\to \infty}(\rho_1(a_n,\mu)+\rho_2(b_n,\nu))=0,
		$$
		i.e.,
		$\displaystyle \lim_{n\to \infty}p_n=(\mu,\nu)$.
	\end{proof}
	
	We will show that some known results are particular case of the main results.
	\begin{ex}\label{e2}
		Theorem \ref{t-1} is a particular case of Theorem \ref{t1}.
	\end{ex}
	\begin{proof}{\rm
			Let $\{A_i\}_{i=1}^3$ be closed, convex subsets of a uniformly convex Banach space $(X, \|\cdot\|)$,  $\displaystyle T :\bigcup_{i=1}^3 A_i\to\bigcup_{i=1}^3 A_i$ be a 3-cyclic summing contraction.
			
			Let $A=A_1\times A_1$, $B=A_2\times A_3$, $C=(A_2\times A_3)\cup\{1\}$, and 
			$$P=\{(g,(b,c),1,(b,c))\in G\times B\times C^2:g\in G,b\in A_2, c\in A_3\},$$
			where $G=\{(a,a)\in A:a\in A_1\}$. Then, $A,B\subseteq (X^2,d)$, where\\
			$d((x_1,x_2),(x_3,x_4))=\|x_1-x_3\|+\|x_2-x_4\|$
			
			Let $T_A:A\times C\to A$ as $T_A((a_1,a_2),p)=(T^3a_1,T^3a_2)$, where $a_1,a_2\in A_1$ and $p\in C$,\\
			$H_A:A\times C\to \{1\}\subseteq C$,\\
			$T_B:B\times C\to B$ as $T_B((b,c),p)=(T^3b,T^3c)$, where $b\in A_2$, $c\in A_3$ and $p\in C$,\\
			$H_B\equiv T_B$,\\
			$f_A:C\to \{0\}\subset \mathbb{R}$, and\\
			$f_B:C\to \mathbb{R}$ as $f_B(b,c)=\|b-c\|$, $f_B(1)= \dist(A_2,A_3)$, where $b\in A_2$ and $c\in A_3$.
			
			It is easy to check $ \dist(A,B)= \dist(A_1,A_2)+ \dist(A_1,A_3)$ and\\
			$\displaystyle \inf_{x\in C}f_B(x)= \dist(A_2,A_3)$.\\
			Thus, for every $(x,y,u,v)=((a,a),(b,c),1,(b,c))\in P$, there holds
			\begin{equation}\label{e2e1}
				\begin{array}{rl}
					&d(T_A(x,u),T_B(y,v))+f_A(H_A(x,u))+f_B(H_B(y,v))\\
					\leq&\lambda(d(x,y)+f_A(u)+f_B(v))+(1-\lambda)S_{f_A,f_B}^C(A,B)\\
					&\Leftrightarrow\\
					&d(T_A((a,a),1),T_B((b,c),(b,c)))+f_B(H_B((b,c),(b,c)))\\
					\leq&\lambda(d((a,a),(b,c))+f_B(b,c))+(1-\lambda)D\\
					&\Leftrightarrow\\
					&d((T^3a,T^3a),(T^3b,T^3c))+f_B(T^3b,T^3c)\\
					\leq&\lambda(\|a-b\|+\|a-c\|+\|b-c\|)+(1-\lambda) D\\
					&\Leftrightarrow\\
					&\|T^3a-T^3b\|+\|T^3a-T^3c\|+\|T^3b-T^3c\|\\
					\leq&\lambda(\|a-b\|+\|a-c\|+\|b-c\|)+(1-\lambda)D,
				\end{array}
			\end{equation}
			where $D=\dist(A_1,A_2)+\dist(A_2,A_3)+\dist(A_3,A_1)$.
			
			Applying (\ref{bz}) $3$ times consecutively for an arbitrary $a\in A_1$, $b\in A_2$, and $c\in A_3$, we get
			$$
			\begin{array}{ll}
				&\|T^3a-T^3b\|+\|T^3a-T^3c\|+\|T^3b-T^3c\|\\
				\leq&k(\|T^2a-T^2b\|+\|T^2a-T^2c\|+\|T^2b-T^2c\|)+(1-k)D\\
				\leq&k^2(\|Ta-Tb\|+\|Ta-Tc\|+\|Tb-Tc\|)+(1-k^2)D\\
				\leq&k^3(\|a-b\|+\|a-c\|+\|b-c\|)+(1-k^3)D.
			\end{array}
			$$
			By the last chain of inequalities and (\ref{e2e1}), we can observe inequality (\ref{bz}) generates inequality (\ref{d1e1}), with $\lambda=k^3$. 
			
			It is easy to check that the definition \ref{contraction}:\ref{C-1} holds. Then, using the fact that inequality (\ref{bz}) generates inequality (\ref{d1e1}), the assumption that $T$ is a 3-cyclic summing contraction, $\displaystyle \inf_{x\in C}f_B(x)= \dist(A_2,A_3)$, and $\displaystyle \inf_{x\in C}f_A(x)=0$, it follows $(T_A,T_B,H_A,H_B,f_A,f_B)\in CEF_C^P(A,B)$.
			
			From example \ref{e0}, we get that the ordered pairs $(A_1,A_2)$ and $(A_1,A_3)$ satisfy the property $UC$. Then, by theorem \ref{t0} and proposition \ref{t0.1}, it follows that the ordered pair $(A,B)$ satisfies the property $CD$. Then, using $(T_A,T_B,H_A,H_B,f_A,f_B)\in CEF_C^P(A,B)$, we can apply theorem \ref{t1}.
			
			\begin{enumerate}[label= con.-\arabic*]
				\item\label{con1}
				For any $g_1\in G$, $g_2\in G$, and $y_0\in B$, the inclusions
				$(g_1,y_0,1,y_0)\in P$ and $(g_2,y_0,1,y_0)\in P$ are true. Then, from theorem \ref{t1} consequence \ref{t1-cons1}, we obtain that $\ds\lim_{n\to\infty}x_n=\alpha \in A$, where $\{x_n\}_{n=0}^{\infty}$ and $\{1\}_{n=0}^{\infty}$ are the iterated sequences of $(T_A,H_A)$ for an arbitrary initial guess $(g,1)\in G\times C$. By theorem \ref{t1} consequence \ref{t1-cons2} it follows $\alpha$ does not depend on the choice of $g$.
				\item\label{con2} 
				Let $\{x_n\}_{n=0}^{\infty}$, $\{1\}_{n=0}^{\infty}$ are the iterated sequences of $(T_A,H_A)$ for an arbitrary initial guess $(g,1)\in G\times C$, and $\ds z=\lim_{n\to\infty}x_n$.
				It is easy to check that $T_A(G\times C)\subseteq G$. Therefore
				\begin{equation}\label{e2e2}
					\{x_n\}_{n=0}^{\infty}\subseteq G.
				\end{equation}
				The metric spaces $(G,d)$ and $(A_1,2\rho)$ are isometric, $(A_1,\rho)$ is completed by assumption. Consequently, $(G,d)$ is also completed. Then, using (\ref{e2e2}), we get $z\in G$.
				For every $g\in G$ and $(y,y)\in B\times C$, there holds $\{1\}_{n=1}^{\infty}\in IS_{f_A}^P(g,y,y)$. Then, from $z\in G$ and theorem \ref{t1} consequences \ref{t1-cons3} and \ref{t1-cons4}, we can observe
				$$\{z\}=\left\{p\in G: \lim_{n\to \infty}T_A(p,1)=p\right\}=\left\{(a,a): a\in A_1,\ T^3a=a\right\}.$$
				\item\label{con3}
				For any $a\in A_1$, $b\in A_2$, and $c\in A_3$, the inclusion\\
				$((a,a),(b,c),1,(b,c))\in P$, is true,\\ 
				$\{(T^{3n}b,T^{3n}c)\}_{n=0}^\infty$, $\{(T^{3n}b,T^{3n}c)\}_{n=0}^\infty$ are the iterated sequences of\\
				$(T_B,H_B)$ for an arbitrary initial guess\\
				$((b,c),(b,c))\in B\times C$, $\displaystyle \inf_{x\in C}f_B(x)= \dist(A_2,A_3)$. Then using theorem \ref{t1} consequences \ref{t1-cons1}, we obtain that for every $b\in A_2$ and $c\in A_3$ there holds $\displaystyle \lim_{n\to \infty}\rho(T^{3n}b,T^{3n}c)=\lim_{n\to \infty}f_B((T^{3n}b,T^{3n}c))= \dist(A_2,A_3)$.
			\end{enumerate}
			
			From \ref{con1} and the definition of $T_A$, it follows that there is $z_1\in A_1$ so that for any arbitrary chosen $a\in A_1$, the sequence $\{T^{3n}a\}_{n=0}^\infty$ converges to $z_1$. By similar arguments it follow that there is $z_i\in A_i$, $i=2,3$ so that for any arbitrary chosen $a_i\in A_i$, $i=2,3$, the sequences $\{T^{3n}a_i\}_{n=0}^\infty$ converge to $z_i$, $i=2,3$, respectively. Then, by the continuity of the metric function $\rho(\cdot,\cdot)$ and \ref{con3}, we obtaine $\rho(z_2,z_3)= dist(A_2,A_3)$. by interchanging the the sets $A_1$, $A_2$ and $A_3$ in the definitions of $A$, $B$, $C$ and $f_B$ we can prove in a similar fashion that $\rho(z_1,z_2 )= dist(A_1,A_2)$ and $\rho(z_1,z_3)= dist(A_1,A_3)$.
			
			Using \ref{con2} and the definition of $T_A$ it follow that $z_1$ is unique fixed point for $T^3$ in $A_1$.
			we can get by similar arguments that $z_i$, $i=2,3$ are the unique fixed points of $T^3$ in $A_i$, for $i=2,3$, respectively. From $T^3(Tz_1)=T(T^3z_1)=Tz_1$ we conclude that $Tz_1=z_2$ and $Tz_2=z_3$ and $Tz_3=z_1$.}
	\end{proof}
	\bibliographystyle{plain}
	\bibliography{bibliography}

\begin{thebibliography}{10}

\bibitem{P-P1}
A.~Abkar and M.~Gabeleh.
\newblock Global optimal solutions of noncyclic mappings in metric spaces.
\newblock {\em Journal of Optimization Theory and Applications},
  153(2):298--305, May 2012.

\bibitem{P-CYCLIC-1}
Hind Alamri, Nawab Hussain, and Ishak Altun.
\newblock Proximity point results for generalized p-cyclic reich contractions:
  An application to solving integral equations.
\newblock {\em Mathematics}, 11(23), 2023.

\bibitem{P-CYCLIC-2}
Victory Asem, Yumnam~Mahendra Singh, Mohammad~Saeed Khan, and Salvatore Sessa.
\newblock On (alpha;,p)-cyclic contractions and related fixed point theorems.
\newblock {\em Symmetry}, 15(10), 2023.

\bibitem{BANACH}
S.~Banach.
\newblock Sur les op´erations dan les ensembles abstraits et leurs
  applications aux integrales.
\newblock {\em Fund. Math.}, 3:133–181, 1922.

\bibitem{CHATTERJEA}
S.~K. Chatterjea.
\newblock Fixed-point theorems.
\newblock {\em C. R. Acad. Bulg. Sci.}, 25:727--730, 1972.

\bibitem{PORDEREDMS}
Binayak~S. Choudhury and N.~Metiya.
\newblock Multivalued and singlevalued fixed point results in partially ordered
  metric spaces.
\newblock {\em Arab Journal of Mathematical Sciences}, 17(2):135--151, 2011.

\bibitem{B-MS-INTRODUCED}
Stefan Czerwik.
\newblock Contraction mappings in b-metric spaces.
\newblock {\em Acta Mathematica et Informatica Universitatis Ostraviensis},
  1(1):5--11, 1993.

\bibitem{Guo}
V.~Lakshmikantham D.~Guo.
\newblock Coupled fixed points of nonlinear operators with applications.
\newblock {\em Nonlinear Analysis, Theory, Methods and Applications},
  11:623--632, 1987.

\bibitem{NONCYCLIC}
A.~Anthony Eldred, W.~A. Kirk, and P.~Veeramani.
\newblock Proximal normal structure and relatively nonexpansive mappings.
\newblock {\em Studia Mathematica}, 171:283--293, 2005.

\bibitem{BPP}
A.~Anthony Eldred and P.~Veeramani.
\newblock Existence and convergence of best proximity points.
\newblock {\em Journal of Mathematical Analysis and Applications},
  323(2):1001--1006, 2006.

\bibitem{UC-5}
Nilakshi Goswami and Raju Roy.
\newblock Some coupled best proximity point results for weak gkt cyclic
  $\phi$-contraction mappings on metric spaces.
\newblock {\em Proceedings of the Jangjeon Mathematical Society},
  23(4):485--502, 2021.

\bibitem{COUPLEDFP}
Dajun Guo and V.~Lakshmikantham.
\newblock Coupled fixed points of nonlinear operators with applications.
\newblock {\em Nonlinear Analysis: Theory, Methods Applications},
  11(5):623--632, 1987.

\bibitem{HARDY-ROGERS}
G.~E. Hardy and T.~D. Rogers.
\newblock A generalization of a fixed point theorem of reich.
\newblock {\em Canadian Mathematical Bulletin}, 16(2):201–206, 1973.

\bibitem{UCBS1}
Miroslav Hristov, Atanas Ilchev, Petar Kopanov, Vasil Zhelinski, and Boyan
  Zlatanov.
\newblock Best proximity points for p ndash;cyclic infimum summing
  contractions.
\newblock {\em Axioms}, 12(7), 2023.

\bibitem{P-CYCLIC-4}
Miroslav Hristov, Atanas Ilchev, Diana Nedelcheva, and Boyan Zlatanov.
\newblock Existence of coupled best proximity points of p-cyclic contractions.
\newblock {\em Axioms}, 10(1), 2021.

\bibitem{ERROR-1}
A.~Ilchev.
\newblock On an application of coupled best proximity points theorems for
  solving systems of linear equations.
\newblock In {\em AIP Conference Proceedings}, volume 2048, page 050003, 2018.

\bibitem{ERROR-2}
A.~Ilchev and B.~Zlatanov.
\newblock Error estimates for approximation of coupled best proximity points
  for cyclic contractive maps.
\newblock {\em Applied Mathematics and Computation}, 290:412--425, 2016.

\bibitem{Ilchev20172873}
Atanas Ilchev and Boyan Zlatanov.
\newblock Fixed and best proximity points for kannan cyclic contractions in
  modular function spaces.
\newblock {\em Journal of Fixed Point Theory and Applications}, 19(4):2873 –
  2893, 2017.
\newblock Cited by: 12.

\bibitem{KANAN}
R.~Kannan.
\newblock Some results on fixed points—ii.
\newblock {\em The American Mathematical Monthly}, 76(4):405--408, 1969.

\bibitem{PMS-3}
Erdal Karapinar, Ibrahim~M Erhan, and A~Yildiz~Ulus.
\newblock Fixed point theorem for cyclic maps on partial metric spaces.
\newblock {\em Applied Mathematics and Information Sciences}, 6(2):239--244,
  May 2012.

\bibitem{PCYCLIC}
S.~Karpagam and Sushama Agrawal.
\newblock Best proximity point theorems for p-cyclic meir-keeler contractions.
\newblock {\em Fixed Point Theory and Applications}, 2009(1):197308, Jan 2009.

\bibitem{P-CYCLIC-5}
Saravanan Karpagam.
\newblock Existence of fixed points and best proximity points of p-cyclic
  boyd-wong contractions.
\newblock {\em Annals of the Academy of Romanian Scientists: Series on
  Mathematics and its Applications}, 13(1-2):111--121, 2021.

\bibitem{Khamsi-Kozlowski}
Kozlowski~W.M. Khamsi, M.A.
\newblock {\em Fixed point theory in modular function spaces}.
\newblock 2015.

\bibitem{KKR}
Kozlowski W.M. Reich~S. Khamsi, M.A.
\newblock Fixed point theory in modular function spaces.
\newblock {\em Nonlinear Analysis}, 14:935--953, 1990.

\bibitem{CYCLIC}
W.~Kirk, P.~Srinivasan, and P.~Veeramani.
\newblock Fixed points for mappings satisfying cyclical contractive conditions.
\newblock {\em Fixed Point Theory}, 4:79--189, 2003.

\bibitem{Kozlowski}
W.M. Kozlowski.
\newblock Advancements in fixed point theory in modular function spaces.
\newblock {\em Arabian Journal of Mathematics}, 1:477--494, 2012.

\bibitem{UC-4}
Ashish Kumar.
\newblock {\em Metric fixed point theory in context of cyclic contractions},
  pages 151--216.
\newblock June 2021.

\bibitem{PMS-2}
Santosh Kumar and Johnson~Allen Kessy.
\newblock Fixed point theorems for hybrid pair of weak compatible mappings in
  partial metric spaces.
\newblock {\em Mathematica Bohemica}, 148(2):223--236, 2023.

\bibitem{ZLATANOV}
B.~Zlatanov M.~Petric.
\newblock Best proximity points and fixed points for p-summing maps.
\newblock {\em Fixed Point Theory and Applications}, 2012:1687--1812, 2012.

\bibitem{P-CYCLIC-3}
P.~Magadevan, S.~Karpagam, and E.~Karapınar.
\newblock Existence of fixed point and best proximity point of p-cyclic orbital
  phi-contraction map.
\newblock {\em Nonlinear Analysis: Modelling and Control}, 27(1):91--101, 2022.

\bibitem{PMS-INTRODUCTED}
S.~G. Matthews.
\newblock Partial metric topology.
\newblock {\em Annals of the New York Academy of Sciences}, 728(1):183--197, 11
  1994.

\bibitem{MK}
A.~Meir and Emmett Keeler.
\newblock A theorem on contraction mappings.
\newblock {\em Journal of Mathematical Analysis and Applications},
  28(2):326--329, 1969.

\bibitem{NONCYCLIC3}
Edraoui Mohamed, Aamri Mohamed, and Lazaiz Samih.
\newblock Relatively cyclic and noncyclic p-contractions in locally k-convex
  space.
\newblock {\em Axioms}, 8(3), 2019.

\bibitem{B-MS-1}
Sushanta~Kumar Mohanta and Ratul Kar.
\newblock Some fixed point results in an ordered b-metric space with an
  application to nonlinear matrix equation.
\newblock {\em Indian Journal of Mathematics}, 65(2):199--227, August 2023.

\bibitem{PMS-1}
Wudthichai Onsod, Poom Kumam, and Teerapol Saleewong.
\newblock Fixed point of suzuki-geraghty type o-contractions in partial metric
  spaces with some applications.
\newblock {\em Science and Technology Asia}, 28(1):17--25, January-March 2023.

\bibitem{ZLATANOV2}
Mihaela Petric and Boyan Zlatanov.
\newblock Best proximity points for p-cyclic summing iterated contractions.
\newblock {\em FILOMAT}, 32(9):3275--3287, 2018.

\bibitem{PETROSEL}
Adrian Petru{\c{s}}el.
\newblock Fixed points vs. coupled fixed points.
\newblock {\em Journal of Fixed Point Theory and Applications}, 20(4):150, Oct
  2018.

\bibitem{B-MS-2}
Adrian Petru{\c{s}}el and Gabriela Petru{\c{s}}el.
\newblock Graphical contractions and common fixed points in b-metric spaces.
\newblock {\em Arabian Journal of Mathematics}, 12(2):423--430, Aug 2023.

\bibitem{Koleva}
B.~Zlatanov R.~Koleva.
\newblock On fixed points for chatterjea’s maps in b-metric spaces.
\newblock {\em Turkish Journal of Analysis and Number Theory}, 4:31--34, 2016.

\bibitem{Ran}
Reurings~MCB Ran, ACM.
\newblock A fixed point theorem in partially ordered sets and some applications
  to matrix equations.
\newblock {\em Proc. Amer. Math. Soc.}, 132:1435--1443, 2004.

\bibitem{NONCYCLIC1}
A.~Safari-Hafshejani.
\newblock Error estimates for approximating fixed points and best proximity
  points for noncyclic and cyclic contraction mappings.
\newblock {\em Iranian Journal of Numerical Analysis and Optimization},
  13(3):385--396, 2023.

\bibitem{UC-3}
Akram Safari-Hafshejani.
\newblock The existence of best proximity points for generalized cyclic
  quasi-contractions in metric spaces with the uc and ultrametric properties.
\newblock {\em Fixed Point Theory}, 23(2):507--518, June 2022.

\bibitem{NONCYCLIC2}
Akram Safari-Hafshejani.
\newblock Optimal common fixed point results in complete metric spaces with
  w-distance.
\newblock {\em Sahand Communications in Mathematical Analysis}, 19(4):117--132,
  2022.

\bibitem{P-P2}
V.~{Sankar Raj}.
\newblock A best proximity point theorem for weakly contractive
  non-self-mappings.
\newblock {\em Nonlinear Analysis: Theory, Methods and Applications},
  74(14):4804--4808, 2011.

\bibitem{B-MS-3}
Syed Shah~Khayyam, Muhammad Sarwar, Asad Khan, Nabil Mlaiki, and Fatima~M.
  Azmi.
\newblock Solving integral equations via fixed point results involving
  rational-type inequalities.
\newblock {\em Axioms}, 12(7), 2023.

\bibitem{UCSTAR}
Wutiphol Sintunavarat and Poom Kumam.
\newblock Coupled best proximity point theorem in metric spaces.
\newblock {\em Fixed Point Theory and Applications}, 2012(1):93, Jun 2012.

\bibitem{UC-2}
Thounaojam Stephen, Rohen Yumnam, Hüseyin Işık, and Laishram Shanjit.
\newblock Some best proximity point results for generalized cyclic contraction
  mappings.
\newblock 03 2023.

\bibitem{UC}
Tomonari Suzuki, Misako Kikkawa, and Calogero Vetro.
\newblock The existence of best proximity points in metric spaces with the
  property uc.
\newblock {\em Nonlinear Analysis: Theory, Methods and Applications},
  71(7):2918--2926, 2009.

\bibitem{Lakshmikantham}
V.~Lakshmikantham T.~Gnana~Bhaskar.
\newblock Fixed point theorems in partially ordered metric spaces and
  aplications.
\newblock {\em Nonlinear Analysis}, 2006:1379--1393, 65.

\bibitem{NONCYCLIC4}
Kazimierz Włodarczyk, Robert Plebaniak, and Artur Banach.
\newblock Best proximity points for cyclic and noncyclic set-valued relatively
  quasi-asymptotic contractions in uniform spaces.
\newblock {\em Nonlinear Analysis: Theory, Methods abd Applications},
  70(9):3332--3341, 2009.

\bibitem{NONCYCLIC5}
Kazimierz Włodarczyk, Robert Plebaniak, and Cezary Obczyński.
\newblock Convergence theorems, best approximation and best proximity for
  set-valued dynamic systems of relatively quasi-asymptotic contractions in
  cone uniform spaces.
\newblock {\em Nonlinear Analysis: Theory, Methods and Applications},
  72(2):794--805, 2010.

\bibitem{UC-UCZ}
V.~Zhelinski and B.~Zlatanov.
\newblock On the uc and uc star properties and the existence of best proximity
  points in metric spaces.
\newblock {\em Annuaire de l’Université de Sofia “St. Kliment Ohridski”.
  Faculté de Mathématiques et Informatique}, 109:121--146, 2022.

\bibitem{UC-1}
Mi~Zhou, Nicolae~Adrian Secelean, Naeem Saleem, and Mujahid Abbas.
\newblock Best proximity points for alternative p-contractions.
\newblock {\em Journal of Inequalities and Applications}, 2024(1):4, Jan 2024.

\bibitem{ZLATANOV-m}
B.~Zlatanov.
\newblock Best proximity points in modular function spaces.
\newblock {\em Arabian Journal of Mathematics}, 4:215--227, 2015.

\bibitem{ERROR-3}
B.~Zlatanov.
\newblock Coupled best proximity points for cyclic contractive maps and their
  applications.
\newblock {\em Fixed Point Theory}, 22(1):431--452, 2021.
\newblock Web of Science, IF=1.287, Q1; SCOPUS, SJR=0.68, Q2, MR4269039, Zbl
  07370686.

\end{thebibliography}
\end{document}